\documentclass[preprint,twocolumn,5p,authoryear]{elsarticle}
  
  \journal{Automatica}
\usepackage{graphicx}
\usepackage{amsmath,amssymb,amsfonts, amsthm}
\usepackage{bbm} 
\usepackage{tikz}
\usetikzlibrary{calc}
\usepackage{subfigure}
\usepackage{verbatim}

\usepackage[boxruled]{algorithm2e}

\usepackage[font=footnotesize]{caption}

\usepackage{sublabel}
\usepackage{enumitem}
\newlist{hyp}{enumerate}{1}
\setlist[hyp]{label = {\bf(H\arabic*)}, resume}
\newlist{ass}{enumerate}{1}
\setlist[ass]{label = {\bf(A\arabic*)}, resume}
\newlist{sol}{enumerate}{1}
\setlist[sol]{label = {\bf(S\arabic*)}, resume}
\newlist{eprop}{enumerate}{1}
\setlist[eprop]{leftmargin=+4em, label = {\bf(EP-\arabic*)}, resume}

\newtheorem{Theorem}{Theorem}
\newtheorem{Definition}[Theorem]{Definition}
\newtheorem{Proposition}[Theorem]{Proposition}
\newtheorem{Lemma}[Theorem]{Lemma}
\newenvironment{Proof}{\textbf{Proof.}}{~\hfill\qed}
\newtheorem{Remark}[Theorem]{Remark}
\newtheorem{Example}[Theorem]{Example}
\newtheorem{Assumption}{Assumption}

\newcommand{\me}{\mathtt{e}}

\newcommand{\D}{\mathbb{D}}

\newcommand{\N}{\mathbb{N}}

\newcommand{\R}{\mathbb{R}}

\newcommand{\cC}{\mathcal{C}}

\newcommand{\cF}{\mathcal{F}}

\newcommand{\cI}{\mathcal{I}}

\newcommand{\cK}{\mathcal{K}}
\newcommand{\cL}{\mathcal{L}}
\newcommand{\cM}{\mathcal{M}}
\newcommand{\cN}{\mathcal{N}}
\newcommand{\cO}{\mathcal{O}}

\newcommand{\cW}{\mathcal{W}}

\newcommand{\fO}{\mathfrak{O}}

\DeclareMathOperator{\rk}{rank}
\DeclareMathOperator{\im}{im}

\DeclareMathOperator{\linspan}{span}

\newcommand{\setdef}[2]{\left\{#1\,\left|\,\vphantom{#1} #2\,\right.\!\!\right\}}

\newcommand{\pw}{\text{\normalfont pw}}

\newcommand{\imp}{\text{\normalfont imp}}
\newcommand{\obs}{\text{\normalfont obs}}
\newcommand{\diff}{\text{\normalfont diff}}
\newcommand{\myleft}{\text{\normalfont left}}

\newcommand{\new}{\text{\normalfont new}}
\newcommand{\ideal}{\text{\normalfont ideal}}

\newcommand{\eps}{\varepsilon}

\newcommand{\Dpwsm}{\D_{\pw\cC^\infty}}

\newcommand{\du}{\mathtt{u}}
\newcommand{\dy}{\mathtt{y}}

\newenvironment{smallpmatrix}
{\left(\begin{smallmatrix}}
{\end{smallmatrix}\right)}

\newenvironment{smallbmatrix}
{\left[\begin{smallmatrix}}
{\end{smallmatrix}\right]}

\begin{document}

\begin{frontmatter}

\title{Detectability and Observer Design for Switched Differential Algebraic Equations\tnoteref{dfg}}
\tnotetext[dfg]{This work was supported by DFG-project TR 1223/2-1 and was partly carried out while the authors were at the University of Kaiserslautern, Germany. The first author also acknowledges the support provided by the ANR project {\sc ConVan} with grant number ANR-17-CE40-0019-01.}
\author[A]{Aneel Tanwani}\ead{aneel.tanwani@laas.fr }
\author[S]{Stephan Trenn}\ead{s.trenn@rug.nl }

\address[A]{Team MAC (Decision and Optimization), LAAS -- CNRS, Toulouse, France}
\address[S]{Jan C.\ Willems Center for Systems and Control, University of Groningen, Netherlands.}   



\begin{abstract}
This paper studies detectability for switched linear differential-algebraic equations (DAEs) and its application to the synthesis of observers, which generate asymptotically converging state estimates.
Equating detectability to asymptotic stability of zero-output-constrained state trajectories, and building on our work on interval-wise observability, we propose the notion of interval-wise detectability: If the output of the system is constrained to be identically zero over an interval, then the norm of the corresponding state trajectories scales down by a certain factor at the end of that interval. Conditions are provided under which the interval-wise detectability leads to asymptotic stability of zero-output-constrained state trajectories.
An application is demonstrated in designing state estimators. Decomposing the state into observable and unobservable components, we show that if the observable component of the system is reset appropriately and persistently, then the estimation error converges to zero asymptotically under the interval-wise detectability assumption.
\end{abstract}

\end{frontmatter}

\section{Introduction}

The growing application of switched systems in modeling and analysis has contributed toward immense research in the area of dynamical systems which combine discrete and continuous dynamics. Different classes of switched systems can be introduced based on the models associated with the switching signal, or the particular characteristics of the individual subsystems. In this regard, this article studies the problem of {\em detectability} for switched systems where the subsystems are described by differential-algebraic equations (DAEs) and the switching signal is assumed to be known a priori.

Switched DAEs arise naturally when the system dynamics undergo sudden structural changes (switches) and the dynamics of each mode are algebraically constrained \citep{Tren12}. A typical example are electrical circuits with switches where the constraints are induced by Kirchhoff's laws.
Our previous works on structural properties of switched DAEs has addressed the problem of observability \citep{TanwTren12} and the observer design \citep{TanwTren13,TanwTren17a} under the stronger assumption of determinability (which in the nonswitched case is equivalent to observability and roughly speaking means that the state at the end of the observation interval can be determined to any given accuracy). Building on this line of work, this article proposes the notion of (interval-) detectability for switched linear DAEs and its application in the observer design.

Roughly speaking, the property of detectability incorporates the notions of observability and stability, that is, a dynamical system is called detectable if the state trajectories, which correspond to the same input and output, converge asymptotically towards a single trajectory. Seen as a generalization of the observability property for classical linear systems, detectability is characterized by asymptotic stability of the unobservable modes for linear time-invariant systems, or stability of the reduced-order system obtained by setting the output of the system to identically zero. For nonlinear systems, while there are different notions for observability \citep{Sont98a}, the notion of output-to-state stability (OSS) provides one possible framework \citep{SontWang97} to study detectability, which has also been used in observer design \citep{AstoPral03}. These techniques are generalized for switched systems as well: The work of \citet{DeSaDiBe09} proposes detectability conditions for switched linear system in terms of the stability of a reduced order switched system. More recently, \citet{MancGarc18} also show the relevance of detectability of a reduced-order system with zero output in establishing global asymptotic stability of the switched system. The notion of OSS has been studied for switched nonlinear systems by \citet{MullLibe12}, where the focus is on characterizing a class of switching signals under which the growth of the state trajectory is bounded by some increasing function of the output norm.

One major utility of the detectability notion is its application in design of observers, or state estimators.\footnote{While some references differentiate between the terms {\em observer}, {\em asymptotic observer} and {\em state estimator}, e.g. \citep{TrenStoo01}, these terms are used synonymously in this article. See the beginning of Section~\ref{sec:observerDesign} for a formal definition adopted in this paper.} The observer design for (nonswitched) DAEs using observability and detectability notions is an ongoing research topic \citep{BergReis17c}. Our approach towards observer design for switched DAEs builds on the observability notions studied in \citep{TanwShim13} and \citep{TanwTren12, TanwTren17a}, where we use the output information from different modes active over an interval to recover the value of the state, either at the start of the interval (observability), or at the end of the interval (determinability). Due to this generalized notion, we have to introduce unobservable dynamics over an interval which not only depend on the unobservable dynamics of individual subsystems but also their activation times. The detectability notion proposed in this article thus relates to the stability of the unobservable dynamics over an interval (and not the individual subsystems). Inspired by the fact that detectability is a sufficient condition for designing state observers for linear systems, we use these ideas to propose an observer design for switched DAEs.

The contribution of this paper lies in studying detectability notions for switched DAEs (see Section~\ref{sec:notionDetect} for the formal definitions) and design state estimators for systems satisfying the detectability assumption in an appropriate sense. This work builds on our two conference papers: geometric conditions for detectability of switched DAEs were studied in \citep{TanwTren15} and the preliminary design of the observer was proposed in \citep{TanwTren17b}. Using the presentation of the later article as a template, this paper provides additional details, rigorous proofs of the results, and simulation results which were not a part of the conference paper. To the best of our knowledge, these results are also new for the case of switched ordinary differential equations (switched ODEs), as the previous works have only dealt with observable switched systems \citep{TanwShim13}. It turns out that an observer for the detectable case has to work fundamentally different to our observer proposed for the determinable case. We illustrate this by the following simple example.

\begin{Example}\label{ex:swODEdim3}
Consider the switched ODE on the interval $[0,3)$ given by
\[
\setlength{\arraycolsep}{0.2ex}
\begin{array}{rcl|rcl}
\dot{x}_1(t) &=& 0 &\quad \dot{x}_1(t) &=& x_2(t) \\
\dot{x}_2(t) &=& 0 &\quad \dot{x}_2(t) & = & 0\\
\dot{x}_3(t) &=& 0 &\quad \dot{x}_3(t) &=& x_2(t) - x_3(t)\\[1ex]
y(t) &=& x_1(t) &\quad y(t) &=& 0 \\[2ex]
\multicolumn{3}{l|}{t\in[0,1)\cup[2,3), \quad} & \multicolumn{3}{l}{\quad t\in[1,2).}
\end{array}
\]
If we restrict our attention to the interval $[0,3)$, then $y(t) \equiv 0$ on this interval implies $x_1(t) \equiv x_2(t) \equiv 0$, and hence $(x_1,x_2)$ is observable (but only when two switches occur, otherwise $x_2$ is not observable). Also, the identically zero output would imply that the magnitude of $x_3$ decreases, which is the notion of detectability we adopt in this paper (see Section~\ref{sec:notionDetect}). It is possible to design an impulsive estimator with states $\widehat x_1, \widehat x_2, \widehat x_3$ which copies the system dynamics over the interval $[0,3)$, and at $t=3$ we reset the estimations of the observable states as
\[
\begin{pmatrix} \widehat x_1(3) \\ \widehat x_2(3) \end{pmatrix} := \fO (y_{[0,3)})
\]
for some map $\fO$, so that, if $e = \widehat x - x$ denotes the state estimation error, we have
\[
\left\vert \begin{pmatrix} e_1(3) \\ e_2(3) \end{pmatrix} \right\vert
\le \alpha \left\vert \begin{pmatrix} e_1(0) \\ e_2(0) \end{pmatrix} \right\vert
\]
for some desired $\alpha \in (0,1)$. Here, and in the remainder of this article, we use the notation $\vert \cdot \vert$ to denote the Euclidean norm of a vector.
Moreover, for the unobservable error $e_3$, we get
\begin{equation}\label{eq:errx_3}
\begin{aligned}
\dot{e}_3(t) &= 0, &\quad &t \in [0,1) \cup [2,3) \\
\dot{e}_3(t) & = e_2(t) - e_3(t), &\quad& t \in [1,2)
\end{aligned}
\end{equation}
and hence
\[
e_3(3) = \me^{-1}e_3(0) + (1 - \me^{-1})e_2(0).
\]
Thus, independently of the accuracy of the estimation of the observable components, for a large initial value $e_2(0)$, the final error $e_3(3)$  may be significantly larger than $e_3(0)$. Therefore a direct application of our previous presented observer to detectable systems will not work. The underlying problem for this example is that it is not enough to have a good estimate of the observable states at the end of the considered interval, but the estimate must be available already when the observable states influence the unobservable states.

\end{Example}

The remainder of the paper is structured as follows: In Section \ref{sec:systemclass} we formally introduce the system class of switched DAEs and also highlight the importance of taking induced Dirac impulses into account. Afterwards we introduce in Section~\ref{sec:notionDetect} the notion of detectability. In particular, we introduce the notion of uniform interval-detectability, which is fundamental for our observer design, which we present in Section~\ref{sec:observerDesign}. The key result in Section~\ref{sec:observerDesign} is Theorem \ref{thm:xi_correction}, which shows how the ideal correction term decreases the estimation error. Convergence of the observer for non-ideal correction terms is shown in Theorem~\ref{thm:convergence} in Section~\ref{sec:convergence}. The observer design in the form of an algorithm and implementational issues are discussed in Section~\ref{sec:implementation}, simulations are carried in Section~\ref{sec:simulations}.

\section{Preliminaries}\label{sec:systemclass}
\subsection{Switched DAEs}
We consider switched linear DAEs of the form
\begin{equation}\label{eq:swDAE}
\begin{aligned}
E_\sigma \dot x &= A_\sigma x + B_\sigma u \\
y &= C_\sigma x + D_\sigma u
\end{aligned}
\end{equation}
where $x, u , y$ denote the state (with dimension $n\in\N$), input (with dimension $\du\in\N$) and output (with dimension $\dy\in\N$) of the system,  respectively. The switching signal $\sigma: [0,\infty) \rightarrow \N$ is a piecewise constant, right-continuous function of time and in our notation it changes its value at time instants $0<t_1<t_2<\ldots$ called the switching times.
We adopt the convention that over the interval $[t_k, t_{k+1})$ of length $\tau_k:=t_{k+1}-t_k$, the active mode is defined by the quintuple $(E_k,A_k,B_k,C_k,D_k)\in\R^{n\times n} \times \R^{n\times n} \times \R^{n\times \du} \times \R^{\dy\times n}\times \R^{\dy\times \du}$, $k \in \N$ and $t_0:=0$. If $E_k=I$ for all $k\in\N$ we call \eqref{eq:swDAE} a switched ODE. In general, $E_k$ is not assumed to be invertible, which means that in addition to differential equations the state $x$ has to satisfy certain algebraic constraints. At a switching instant the algebraic constraints before the switch and the algebraic constraints after the switch do not match in general, i.e.\ the state variable has to jump in order to satisfy the algebraic constraints after the switch (cf.\ Lemma~\ref{lem:consProj} in Appendix~\ref{sec:basicProp}). These induced jumps are a first major difference to switched ODEs (which do not exhibit jumps unless one imposes some additional external jump rules). The second major difference is the possible presence of Dirac impulses in the state variable $x$ in response to a state jump, see \citet{Tren12} for details. The following example shows this effect and also the role of the Dirac impulses in the state estimation problem.
\begin{Example}\label{ex:impulses_important}
   Consider the switched DAE given by, $i\in\N$,
   \[
   \setlength{\arraycolsep}{0.2ex}
\begin{array}{rcl|rcl}
t&\in&[2i,2i+1) \quad&\quad t&\in&[2i+1,2i+2) \\[2ex]
\dot{x}_1 &=& x_3 & 0 &=& x_1 \\
\dot{x}_2 &=& 0 & \dot{x}_1 &=& x_2\\
\dot{x}_3 &=& 0 & \dot{x}_3 &=& 0 \\
\dot{x}_4 &=& x_3-x_4 & \dot{x}_4 &=& x_3 - x_4\\[1ex]
y &=& 0 &\quad y &=& x_2
\end{array}
   \]
The dynamics for $x_3$ and $x_4$ are actually non-switched and it is therefore obvious that the overall switched system can only be detectable when it is possible to determine $x_3$ from the output. On the intervals $[2i,2i+1)$, $i\in\N$, the output is zero by definition and on the open intervals $(2i+1,2i+2)$, $i\in\N$, it holds that $x_1=0$, hence $y=x_2=\dot{x}_1 =0$, i.e.\ the output is zero almost everywhere. Consequently, we are not able to deduce anything about $x_3$ from the output if we do not take into account what the output is doing \emph{at} the switching times. So what is $x_2$ doing \emph{at} the switching times $t=2i+1$, $i\in\N$? The state $x_1$ jumps from $x_1((2i+1)^-)$ to $x_1((2i+1)^+)=0$, hence $x_2$ contains the derivative of this jump! The derivative of a jump is only well defined in a distributional (generalized functions) solution framework; in this framework $x_2$ contains a Dirac impulse with magnitude $-x_1((2i+1)^-)$ and this Dirac impulse is visible at the output. Consider now the switching time $t=3$ then we can deduce from the Dirac impulse of the output at $t=3$ the value $x_1(3^-)$. We know that $x_1=0$ on $(1,2)$, $\dot{x}_1=x_3$ on $[2,3)$ and $x_3$ is constant, hence $x_1(3^-)=x_1(2^-) + (3-2)x_3 = x_3$. In particular, if we observe a zero output (including zero Dirac impulses) we can conclude that $x_3=0$, $(x_1,x_2)=(0,0)$ on $(1,\infty)$ and $x_4(t)\to 0$ as $t\to\infty$. In summary, the above switched DAE is detectable but it is not possible to estimate the state without taking into account the Dirac impulses in the output. 
\end{Example}

The above example shows that an observer design which does not utilize the information from possible Dirac impulses in the output will not work for general switched DAEs. We will therefore recall now the distributional solution framework for \eqref{eq:swDAE} as introduced in \cite{Tren09d}.

\subsection{Distributional solution framework}
Let $\D$ denote the space of distributions in the sense of \cite{Schw50}, i.e.\ $D\in\D$ if, and only if, $D:\cC^\infty_0\to\R$ is linear and continuous, where $\cC^\infty_0$ is the space of test functions consisting of smooth functions $\varphi:\R\to\R$ with compact support and equipped with a suitable topology. Any locally integrable function $f:\R\to\R$ induces a distribution $f_\D\in\D$ given by
\[
   f_\D(\varphi) := \int_\R f \varphi.
\]
For differentiable $f$ it is easily seen via integration by parts that
\[
   (f')_\D(\varphi) = -f_\D(\varphi'), 
\]
which motivates the definition of the derivative of a general distribution $D\in\D$:
\[
   D'(\varphi) := -D(\varphi).
\]
For some interval $\cI\subseteq \R$ let $\mathbbm{1}_{\cI}$ be the indicator function of $\cI$, i.e.\ $\mathbbm{1}_{\cI}(t)=1$ for $t\in \cI$ and $\mathbbm{1}_{\cI}(t)=0$ otherwise. Then the Dirac impulse $\delta$ can be defined as the distributional derivative of the unit jump (or Heaviside step function), i.e.\
\[
   \delta:=((\mathbbm{1}_{[0,\infty)})_\D)'.
\]
Note that $\delta(\varphi) = \varphi(0)$ for any test function $\varphi\in\cC^\infty_0$. The Dirac impulse at $t\in\R$, denoted by $\delta_t$, is the distributional derivative of $\mathbbm{1}_{[t,\infty)}$. As shown in \citep{Tren09d} it is not possible to use the space $\D$ directly as the underlying solution space for the switched DAE \eqref{eq:swDAE}. Instead, the smaller space of piecewise-smooth distributions $\D_{\pw\cC^\infty}$ will be used as underlying solution space for \eqref{eq:swDAE}, where
\[
   \D_{\pw\cC^\infty}\! :=\! \setdef{D\!=\!f_\D+\!\sum_{t\in T}D_t}{\begin{aligned}
     &f\!\in\!\cC^\infty_\pw,\ T\!\subseteq\!\R\text{ discrete},\\
      &D_t\in\linspan\{\delta_t,\delta_t',\delta_t'',...\}\end{aligned}}\!,
\]
i.e.\ a piecewise-smooth distribution is the sum of a piecewise-smooth function and Dirac impulses (and their derivatives) at isolated points in time. Any piecewise smooth distribution $D=f_\D+\sum_{t\in T} D_t\in\D_{\pw\cC^\infty}$ can be evaluate at time $t\in\R$ in three different ways:
\[\begin{aligned}
    D(t^+) &:= f(t^+) := \lim_{\eps \searrow 0} f(t+\eps),\\
    D(t^-) &:= f(t^-) := \lim_{\eps \searrow 0} f(t-\eps),\\
    D[t] &:= \begin{cases}
       D_t,& t\in T\\
       0,& \text{otherwise}
    \end{cases}
  \end{aligned}
\]
where we denote by $D[t]$ the impulsive part of $D$ at time $t$. Note that for any $t\in T$
\[
   D[t] = \sum_{j=0}^{n_t} \alpha_j^t \delta_t^{(j)}
\]
for some finite $n_t\in\R$ and $\alpha_0^t,\alpha_1^t,\ldots,\alpha_{n_t}^t\in\R$. Furthermore, the product of a piecewise-smooth function with a piecewise-smooth distribution is well defined, in particular, \eqref{eq:swDAE} can be evaluated for piecewise-smooth distributions. It is also possible to define the restriction to intervals for piecewise-smooth distributions,in particular, for any interval $[a,b)\subseteq\R$ we have
\[
   D_{[a,b)} = \mathbbm{1}_{[a,b)} D.
\]
\begin{Lemma}[cf.\ \cite{Tren09d}]
   Consider the switched DAE \eqref{eq:swDAE} and assume that each matrix pair $(E_p,A_p)$ is regular, i.e.\ $\det(sE_p-A_p)$ is not the zero polynomial. Then for every $u\in\D_{\pw\cC^\infty}^\du$, any $x_0\in\R^n$ and any interval $[a,b)\subseteq [0,\infty)$ there exists $x\in\D_{\pw\cC^\infty}^n$ uniquely defined on $[a,b)$ such that $x(a^-)=x_0$ and \eqref{eq:swDAE} holds as an equation of piecewise-smooth distributions restricted to $[a,b)$.
\end{Lemma}

This motivates the following solution definition of \eqref{eq:swDAE}.

\begin{Definition}[Solution of switched DAE]
  A tuple $(x,u,y)$ (or just $x$ when $u$ and $y$ are clear) is called a solution of \eqref{eq:swDAE} on an interval $\cI$ if $x\in\D_{\pw\cC^\infty}^n$, $u\in\D_{\pw\cC^\infty}^\du$, $y\in\D_{\pw\cC^\infty}^\dy$ and \eqref{eq:swDAE} restricted to $\cI$ holds in the distributional sense. If $\cI=[0,\infty)$ we omit ``on the interval $[0,\infty)$'' in the following.
\end{Definition}

\section{Detectability Notions}
\label{sec:notionDetect}
Roughly speaking, in classical literature on nonswitched systems, a dynamical system is called detectable if, for a fixed input and an observed output, the trajectories starting from every pair of indistinguishable initial states converge to a common trajectory asymptotically. This definition can readily be generalized to DAEs (see e.g.\ the notion of behavioral detectability in \citet[Sec.~9]{BergReis17a}) as well as to switched systems (see e.g.\ \citet[Defn.~2.2]{DeSaDiBe09}); the formal definition for switched DAEs is as follows:
\begin{Definition}\label{def:globalDet}
  The switched DAE \eqref{eq:swDAE} is called \emph{detectable} for a given switching signal $\sigma$, if there exists a class $\cK\cL$ function\footnote{A function $\beta:\R_{\ge 0} \times \R_{\ge 0} \to \R_{\ge 0}$ is called a class $\cK\cL$ function, if 1) for each $t \ge 0$, $\beta(\cdot, t)$ is continuous, strictly increasing, with $\beta (0,t) = 0$; 2) for each $r \ge 0$, $\beta(r,\cdot)$ is decreasing and converging to zero as $t\to\infty$.}%
$\beta:\R_{\ge 0} \times \R_{\ge 0} \rightarrow \R_{\ge 0}$ such that, for any two distributional solutions $(x_1,u,y)$, $(x_2,u,y)$ of \eqref{eq:swDAE} we have
\begin{equation}\label{eq:stabCont}
  |x_1(t^+)-x_2(t^+)| \le \beta (|x_1(0^-)-x_2(0^-)|, t), \quad \forall \, t \ge 0.
\end{equation}
\end{Definition}

Because of linearity the definition can be simplified to the case that $u=0$ and $y=0$, in particular, convergence to zero has only to be checked for the homogeneous system and the initial states in\begin{equation}
\cN^\sigma:= \setdef{ x^0 \in \R^n }{\begin{aligned} & (x,u=0,y=0) \text{ solves \eqref{eq:swDAE}} \\ & \wedge \, x(0^-) = x^0 \end{aligned}},
\end{equation}
or in other words, detectability is the same as asymptotic stability of the switched DAE \eqref{eq:swDAE} with $u=0$ and $y=0$.

\begin{Remark}
In contrast to previous works on stability of switched DAEs \citep{LibeTren09,LibeTren12} we do not require impulse-freeness of solutions for asymptotic stability.\ The reason is that the presence of Dirac impulses may actually help to make certain states observable (cf.\ Example~\ref{ex:impulses_important}), hence the exclusion of Dirac impulses may exclude an important class of problems where Dirac impulses are needed for observability (or detectability). It should also be noted that the magnitude of the Dirac impulses is always proportional to the state value prior to the time the Dirac impulse occurs (cf.\ the explicit expression \eqref{eq:impSol} in the Appendix), i.e.\ when the state converges to zero as $t\to\infty$ the magnitude of the Dirac impulses also converges to zero (under an additional mild boundedness assumption on $(E_k,A_k)$ as $k\to\infty$).
\end{Remark}

Computation of the set $\cN^\sigma$ in general depends on all switching times and the data of all subsystems.
For certain applications, such as state estimation which we discuss later, it may be desirable to work with system data available on finite intervals only, and in that case, Definition~\ref{def:globalDet} may not be suitable.
To overcome this problem, we consider the system behavior on finite intervals, and introduce the notion of interval-detectability:
\begin{Definition}[Interval-detectability]\label{def:localDet}
The switched DAE \eqref{eq:swDAE} is called $[t_p,t_q)$-detectable for a given switching signal $\sigma$, if there exists a class $\cK\cL$ function $\beta: \R_{\ge 0} \times \R_{\ge 0} \rightarrow \R_{\ge 0}$ with
\begin{subequations}\label{eq:localEstDef}
\begin{equation}\label{eq:localEstDefa}
\beta (r,t_q-t_p) < r, \quad \forall \, r > 0
\end{equation}
and for any local solution $(x,u=0,y=0)$ of \eqref{eq:swDAE} on $[t_p,t_q)$ we have
\begin{equation}\label{eq:localEstDefb}
|x(t^+)| \le \beta (|x(t_p^-)|, t- t_p), \quad \forall \, t \in [t_p, t_q).
\end{equation}
\end{subequations}
\end{Definition}

One should be aware, that a solution on some interval is not always a part of a solution on a larger interval. Consequently, detectability does not always imply interval-detectability: The switched system $0=x$ on $[0,t_1)$ and $\dot{x}=0$ on $[t_1,\infty)$ with zero output is obviously detectable (because zero is the only global solution), but it is not interval-detectable on $[t_1,s)$ for any $s>t_1$ because on $[t_1,s)$ there are nonzero solutions which do not converge towards zero.

Furthermore, we would like to emphasize that the interval $[t_p,t_q)$ in general contains multiple switches, i.e.\ it is \emph{not assumed} that the individual modes of the switched systems are detectable.
We need some uniformity assumption to conclude that interval-detectability on each interval of a partition of $[0,\infty)$ implies detectability:
\begin{Assumption}[Uniform interval-detectability]\label{ass:uniform_localdetec}
  Consider the switched system \eqref{eq:swDAE} with switching signal $\sigma$ and switching times $t_k$, $k\in\N$. Assume that there exists a strictly increasing sequence $(q_i)_{i=0}^{\infty}$ with $q_0>0=:q_{-1}$ such that for $p_i:=q_{i-1}$ the system is $[t_{p_i},t_{q_i})$-detectable with $\cK\cL$-function $\beta_i$ for which additionally it holds that
  \begin{subequations}
\begin{align}
\beta_i(r,t_{q_{i}}-t_{p_i}) & \le \alpha \, r, \quad \forall \, r > 0, \forall \, i \in \N,\label{eq:betai_alpha}\\
\beta_i(r,0) & \le M \, r, \quad \forall \, r > 0, \forall \, i \in \N,\label{eq:betai_M}
\end{align}
\end{subequations}
for some uniform $\alpha\in(0,1)$ and $M\geq 1$.
\end{Assumption}

We can now show the following result:
\begin{Proposition}\label{prop:uniform_local_implies_global}
If the switched system \eqref{eq:swDAE} is uniformly interval-detectable in the sense of Assumption~\ref{ass:uniform_localdetec} then \eqref{eq:swDAE} is detectable.
\end{Proposition}

\begin{Proof}
Let for $i\in\N$
\[
   \widehat{\beta}_i(r,t-t_{p_i}) = M r - (t-t_{p_i})\frac{M r (1-\alpha)}{t_{q_i}-t_{p_i}},
\]
i.e.\ for each $r>0$ the function $\widehat{\beta}_i(r,\cdot)$ is linear on $[t_{p_i},t_{q_i})$ and decreasing from $M r$ towards $\alpha M r$. Now let
\[
   \beta(r,t):=\max\left\{\beta_i(\alpha^i r,t-t_{p_i}), \widehat{\beta}_i(\alpha^i r,t-t_{p_i})\right\},
\]
where $i\in\N$ is such that $t\in[t_{p_i},t_{q_i})$. Clearly, for fixed $t$, $\beta(\cdot,t)$ is continuous and strictly increasing. From \eqref{eq:betai_alpha} and $M\geq 1$ it follows that
\[
   \beta(r,t_{q_i}^-) = \max\left\{\beta_i(\alpha^i r,t_{q_i}-t_{p_i}),M\alpha^{i+1} r\right\} = M\alpha^{i+1} r
\]
and, invoking \eqref{eq:betai_M},
\[
   \beta(r,t_{p_i}) = \max\left\{\beta_i(\alpha^i r,0),M\alpha^i r\right\} = M\alpha^{i} r.
\]
 Because $q_{i}=p_{i+1}$, continuity of $\beta(r,\cdot)$ with fixed $r>0$ follows. Furthermore, on each interval $[t_{p_i},t_{q_i})$ the function $\beta(r,\cdot)$ is strictly decreasing as a maximum of two strictly decreasing functions. Additionally, $\beta(r,t_{p_i}) =  M\alpha^{i} r$ with $\alpha\in(0,1)$ implies that $\beta(r,t)$ converges to zero as $t\to\infty$. So $\beta$ is a $\cK\cL$-function and it remains to be shown that $|x(t^+)|
 \leq \beta(|x(t_0^-)|,t)$ for any solution $(x,u=0,y=0)$ of \eqref{eq:swDAE}. First observe, that by \eqref{eq:betai_alpha} and continuity of $\beta_i$ it follows that
 \[
     |x(t_{p_{i+1}}^-)| = |x(t_{q_i}^-)|\leq \beta_i(|x(t_{p_i}^-)|,t_{q_i}-t_{p_i})\leq \alpha |x(t_{p_i}^-)|,
 \]
 hence $|x(t_{p_i}^-)|\leq \alpha^i |x(0^-)|$. Therefore,
 \begin{multline*}
   |x(t^+)| \leq \beta_i(|x(t_{p_i}^-)|,t-t_{p_i}) \leq \beta_i(\alpha^i|x(0^-)|,t-t_{p_i})\\
    \leq \beta(|x(0^-)|,t).
 \end{multline*}
~\hfill~
\end{Proof}


The uniformity conditions \eqref{eq:betai_alpha} and \eqref{eq:betai_M} are both crucial, see \citep[Example~2]{TanwTren17b} for counterexamples.


\begin{Example}[Example~\ref{ex:swODEdim3} revisited]\label{ex:detODEdim3}
Consider the system in Example~\ref{ex:swODEdim3} with periodic switching where the mode sequence and activation times defined for the interval $[0,3)$ are repeated on the interval $[3i,3i+3)$, $i \in \N$.
It can be verified that the resulting system is uniformly interval-detectable, and hence detectable by Proposition~\ref{prop:uniform_local_implies_global}. To see this, we consider the sequence $q_i = 3i+3$ and let
\[
\beta(r,s) := r\me^{2-s}.
\]
The function $\beta$ satisfies the inequalities \eqref{eq:betai_alpha} and \eqref{eq:betai_M}, with $\alpha = \me^{-1}$ and $M = \me^2$, respectively.
The constraint $y \equiv 0$ yields $x_1 = x_2 \equiv 0$, and it can be verified that
\[
\vert x_3(t) \vert \le \me^{2 - (t-3i)} x_3(3i), \quad t \in [3i, 3i+3).
\]
\end{Example}

\begin{Remark}
   Proposition~\ref{prop:uniform_local_implies_global} can actually be seen as a statement about asymptotic stability of switched systems and when it is possible to conclude asymptotic stability from some stability notion on finite intervals. Furthermore, the statement carries over to the nonlinear case without much change, because in the proof we did not exploit the special (linear) form of the switched system \eqref{eq:swDAE}.
\end{Remark}
%

\section{Observer design}\label{sec:observerDesign}

We now turn our attention to designing observers. By definition, an observer for system~\eqref{eq:swDAE} is an operator $\widehat{\fO}$, either static or dynamic, which for each $(x,u,y)$ satisfying \eqref{eq:swDAE}, generates $\widehat x:=\widehat{\fO}(u,y)$ having the property that
\[
  |\widehat x(t^+)-x(t^+)| \le \beta (|\widehat x(0^-)-x(0^-)|, t), \quad \forall \, t \ge 0
\]
for some class $\cK\cL$ function $\beta$. The observer design presented here is an extension of the algorithm proposed in \citep{TanwTren17a} for the determinable case (in particular, the interval-wise observer design), i.e.\ we propose an impulsive observer which consists of a system copy and a correction term which updates the state of the system copy at the end of the detectability interval.

Taking a bird's eye view, the state estimator --- under the uniform interval-detectability (Assumption~\ref{ass:uniform_localdetec}) with detectability intervals $[t_{p_i},t_{q_i})$, $i\in\N$ --- is given by $\widehat{x}:=\sum_{i\in\N} (\widehat{x}_i)_{[t_{p_i},t_{q_i})}$ with
\begin{equation}\label{eq:system_copies_with_xi_i}
\begin{aligned}
   &\hspace{3ex}\left.\begin{aligned}
      E_\sigma \dot{\widehat{x}}_i &= A_{\sigma} \widehat{x}_i + B_\sigma u,\\
      \widehat{y} &= C_\sigma \widehat{x}_i + D_\sigma u,
   \end{aligned}
   \right\} \text{ on }[t_{p_i},t_{q_i}),\\
   &\widehat{x}_{i+1}(t_{q_i}^-) = \widehat{x}_{i}(t_{q_i}^-)-\xi_i.
   \end{aligned}
\end{equation}
where $\xi_i\in\R^n$ is a state estimation correction obtained from the available data on the interval $[t_{p_i},t_{q_i})$ applied at the end of the corresponding interval. Similar to the technique adopted in \citep{TanwTren17a}, the correction term $\xi_i$ is obtained by collecting the local observability data for each mode. However, these local data is combined in a fundamentally different way compared to \citep{TanwTren17a}, because $\xi_i$ is obtained by composing the local observability data backward in time first and then propagating this forward in time under the error dynamics, cf.\ Example~\ref{ex:swODEdim3}.

In particular, a much more complicated algorithm is needed to obtain the correction term at the end of the interval. In fact, it consists of the three following steps which have to be carried out on each of the detectability intervals $[t_{p_i},t_{q_i})$:
\begin{enumerate}[wide, labelwidth=!, labelindent=0pt]
  \item Collect local observability data for each mode synchronous to the system dynamics from the measured input and output over the interval $[t_{p_i,}t_{q_i})$.
  \item\label{item:propagate_back} Propagate back the collected information to obtain an estimation correction $\xi_i^\myleft$ at the beginning of the detectability interval.
  \item\label{item:propagate_forward} Propagate forward the correction term $\xi_i^\myleft$ to obtain the actual estimation correction $\xi_i$ at the end of the interval.
\end{enumerate}

We will now explain each of the steps in detail, for that we drop the index $i$ and just consider the generic detectability interval $[t_p,t_q)$ for some $q>p\geq 0$. It is helpful to introduce the estimation error $e:=\widehat{x}-x$ (which we don't know, because $x$ is not known) and the corresponding output mismatch $y^e:=\widehat{y}-y$ (which we know). It is easily seen that the error is governed by the following homogeneous switched DAE on $[t_p,t_q)$:
\begin{equation}\label{eq:errorDAE}
\begin{aligned}
   E_{\sigma} \dot{e} = A_{\sigma} e, \qquad
   y^e = C_{\sigma} e
 \end{aligned}
\end{equation}
and the idea of the observer is to estimate the error signal $e$ from the measured output mismatch $y^e$. The estimation $\xi$ of $e(t_q^-)$ will then be used to update the state estimation $\widehat{x}$ at $t_q$ to $\widehat{x}(t_q^-)-\xi$; if $\xi\approx e(t_q^-)$ it then holds by definition that
\[
   \widehat{x}(t_q^-)-\xi\approx x(t_q^-).
\]

\begin{Remark}
A key feature of our observer is the consideration of the homogeneous error dynamics \eqref{eq:errorDAE} not only in the analysis but also in the implementation of our observer. In particular, it is not necessary to store the input and output values over a (possibly long) time interval to carry out Steps  \ref{item:propagate_back} and \ref{item:propagate_forward}, see also Remark~\ref{rem:correction_xi^left}. This approach is only possible because the observer consists of a system copy \emph{without} a continuous update of the state estimation based on an output error injection; instead our observer is an impulsive observer in the sense that only at isolated time points the state estimation is updated discontinuously. Another reason \emph{not to use} continuous updates of the state estimations via output error injection is the problem that the observable subspace is not necessarily aligned with the original coordinates. While for the original Example~\ref{ex:swODEdim3} it would be possible to continuously update the estimation of $x_1$ already on the first interval; this update is not possible if we slightly change the example such that the output on the intervals $[0,1)$ and $[2,3)$ takes the form $y=x_1+x_2$. It is easily seen that also with this output the switched system is detectable; however, now it is unclear how the local (one-dimensional) observability information available on the interval $[0,1)$ can be injected continuously to update the state variable in a meaningful way.
\end{Remark}

\subsection{Collecting local observability data for each mode}\label{sec:zk}

For each mode $k$ with $p\leq k\leq q-1$ consider the local unobservable space:
\begin{equation}\label{eq:Wk}
   \cW_k
   := \setdef{e_0\in\R^n}{\begin{aligned}&e(t_k^-)=e_0, \text{where $(e,y^e=0)$}\\&\text{ solves \eqref{eq:errorDAE} on $[t_k,t_{k+1})$}\end{aligned}}
\end{equation}
Defining $\Pi_k$, $O^\diff_k$ and $O^\imp_k$ in terms of $(E_k,A_k)$ as in the Appendix~\ref{sec:basicProp}, it can be shown (cf.\ \citep{TanwTren13,TanwTren17a}) that
\[
   \cW_k = \Pi_k^{-1}\ker O^\diff_k\ \cap\ \ker O^\imp_k.
\]
Note that in general $\Pi_k$ is not invertible and $\Pi_k^{-1}$ stands for the set-valued preimage.

\begin{Remark}[Different definitions of local unobservable space]
In \citep{TanwTren12,TanwTren13,TanwTren17a} slightly different definitions of the local unobservable spaces are used. The difference is based on the different solution interval; in the previous works this interval was $(t_{k-1},t_{k+1})$ or $(t_{k-1},t_{k}]$, while here the interval is $[t_k,t_{k+1})$. As a consequence $\cW_k$ here only depends on the system's parameters of mode $k$ and not on parameters of two modes. The different definitions are motivated by the overall observability notion studied. In our first work we defined the local unobservable space in such a way that all information \emph{around} a single switching time is utilized, in particular, the local unobservable space for a system with a single switch matched the overall unobservable space. Our later works focused on observer design and determinability (i.e.\ the ability to determine the state value at the \emph{end} of the observation interval), for this reason it made sense to consider as local information the continuous output measurement before the current switching time and the Dirac impulses instantaneously induced at that switching time. Here our focus is on observable components of the state at the \emph{beginning} of the observed interval (cf.\ Example~\ref{ex:swODEdim3}) which motivated the definition \eqref{eq:Wk}.
\end{Remark}

If the output mismatch $y^e$ is nonzero then the value of $e$ in \eqref{eq:errorDAE} prior to the switching time $t_k$ can be decomposed as
\[
   e(t_k^-) = W_k w_k + Z_k z_k,
\]
where $\im W_k=\cW_k$ and $\im Z_k = \cW_k^\bot$ and $W_k,Z_k$ are orthonormal matrices. In particular, $z_k=Z_k^\top e(t_k^-)$ is the observable part of the error $e(t_k^-)$ based on the knowledge on the interval $[t_k,t_{k+1})$. It is possible to write the observable part $z_k$ in terms of $y^e$: 
\begin{equation}\label{eq:zk=cOk}
    z_k = \fO_k(y^e_{[t_k,t_{k+1})})
\end{equation}
with some operator $\fO_k$ which evaluates the impulsive part $y^e[t_k]$ as well as the smooth part $y^e_{(t_k,t_{k+1})}$ (possibly depending on the derivatives of $y^e$).
The construction of this ``ideal'' observability operator $\fO_k$ is provided in Appendix \ref{sec:localObsComp}. One may also refer to \citep[Section~5]{TanwTren17a} for a detailed treatment. In practice, only an approximation $\widehat{\fO}_k$ of $\fO_k$ will be available, this will be discussed in Section~\ref{sec:convergence}.

\subsection{Combining local information backwards in time}
Next we want to combine the observable information $z_p, z_{p+1},\ldots,z_{q-1}$ collected on the interval $[t_p,t_q)$ via \eqref{eq:zk=cOk}, to arrive at an expression for $e(t_p^-)$.
To do so, we first quantify the information that can be extracted from the output over an interval $[t_k,t_q)$ by introducing the subspace
\begin{equation}\label{eq:Nkq}
    \cN_{k}^{q}:= \setdef{e_0\in\R^n}{\begin{aligned}&e(t_k^-)=e_0\text{, where $(e,y^e=0)$}\\ &\text{solves \eqref{eq:errorDAE} on $[t_k,t_q)$}\end{aligned}}
\end{equation}
which can be recursively calculated (backwards in time, i.e.\ for $k=q-1,q-2,\ldots,p$), see \eqref{eq:DAESeqUnobs} in the Appendix.
We then decompose the state estimation error just before the interval $[t_k,t_q)$ accordingly:
\begin{equation}\label{eq:MkNk}
   e(t_{k}^-) = M^q_k \mu_k + N^q_k \nu_k
\end{equation}
for some vectors $\mu_k$ and $\nu_k$ of appropriate dimension; here, $M^q_k$ and $N^q_k$ are the matrices with orthonormal columns such that
\[
   \im N^q_k = \cN_{k}^{q}\quad\text{and}\quad \im M^q_k = (\cN_{k}^{q})^\bot.
\]
As shown in Appendix \ref{sec:intObsComp}, there exists a matrix $\cF^q_k$ given in terms of $M^q_{k+1}$, $N^q_{k+1}$, $(E_k,A_k)$ and the duration time $\tau_k=t_{k+1}-t_k$ such that for $p\leq k \leq q-2$
\[
  \mu_k = \cF^q_k\begin{pmatrix} z_k \\ \mu_{k+1} \end{pmatrix}.
\]
and $\mu_{q-1}=z_{q-1}$. Note that by construction, for all $p\leq k\leq q-1$
\[
   e(t_k^-)-M^q_k \mu_k \in \cN_{k}^q.
\]
Now the ideal estimation error correction is
\begin{equation}\label{eq:cO_p^q-1}
\begin{aligned}
    \xi^\myleft &:= M^q_p \mu_p \\
    &= M^q_p\cF^q_p\!\begin{smallpmatrix} z_p & \\ \!\!{\displaystyle\cF_{p+1}^q}\! &\!\! \begin{smallpmatrix} z_{p+1}& \\ {\displaystyle\cF_{p+2}^q}\! &\!\! \begin{smallpmatrix}\ddots & \\ & & {\displaystyle\cF_{q-2}^q} \begin{smallpmatrix} z_{q-2} \\ z_{q-1} \end{smallpmatrix} \end{smallpmatrix} \end{smallpmatrix}\end{smallpmatrix} \\
    &=: \cO_p^{q-1} \mathbf{z}_p^{q-1},
\end{aligned}
\end{equation}
where $\mathbf{z}_p^{q-1} = (z_p/z_{p+1}/\cdots/z_{q-1})$; here the notation $(\cdots/\cdots/\cdots)$ stands for a vector (or matrix) resulting from stacking all entries over each other. In fact, by construction the following is true:
\[
    \boxed{e(t_p^-) - \xi^\myleft \in \cN_p^q \text{ and }\xi^\myleft\in{\cN_p^q}^\bot,}
\]
i.e.\ we are able to obtain the orthogonal projection of $e(t_p^-)$ onto $\cN_p^q$ without actually knowing $e(t_p^-)$.
\subsection{Propagating correction term forward in time}
For the detectability interval $[t_p,t_q)$, let $\xi^\myleft$ be given as above, then let
\begin{equation}\label{eq:defxi}
   \boxed{\xi:=\Phi_{p}^q \xi^{\myleft},}
\end{equation}
where $\Phi_p^p= I$ and $\Phi_p^{k+1}$, $k=p, p+1,\ldots, q-1$ is recursively given by
\begin{equation}\label{eq:Phi_p^k}
    \Phi_{p}^{k+1} = \me^{A^\diff_{k} \tau_{k}} \Pi_{k} \Phi_p^{k}
\end{equation}
with $\Pi_k$ and $A^\diff_k$ are given as in Definition~\ref{def:proj} in the Appendix. In fact, as a consequence from Lemmas \ref{lem:consProj} and \ref{lem:diffProj} in the Appendix, $\Phi_p^q$ is the transition matrix of the homogeneous error DAE \eqref{eq:errorDAE} from $e(t_p^-)$ to $e(t_q^-)$. We then have the following result:

\begin{Theorem}\label{thm:xi_correction}
   Consider the switched DAE \eqref{eq:swDAE} which is detectable on $[t_p,t_q)$ with corresponding $\cK\cL$-function $\beta$. Let $(\widehat{x},\widehat{y})$ be the solution of the system copy
 \begin{equation}\label{eq:system_copy_on_[tp,tq)}
 \begin{aligned}
     E_\sigma \dot{\widehat{x}} &= A_\sigma \widehat{x} + B_\sigma u, \\
     \widehat{y} &= C_\sigma \widehat{x} + D_\sigma u
  \end{aligned}
 \end{equation}
 on $[t_p,t_q)$.
Based on the output mismatch $y^e=\widehat{y} - y$, let
 \[
   \xi = \Phi_p^q \xi^\myleft = \Phi_p^q \cO_p^{q-1} \mathbf{z}_p^{q-1}
 \]
 where $\Phi_p^q$ is given by \eqref{eq:Phi_p^k}, $\cO_p^{q-1}$ is given by \eqref{eq:cO_p^q-1} and  $\mathbf{z}_p^{q-1} = (z_p/z_{p+1}/\cdots/z_{q-1})$ with $z_k = \fO_k(y^e_{[t_k,t_k+1)})$, $k=p,p+1,\ldots,q-1$ is given by \eqref{eq:zk=cOk}. Then 
 \[
 \begin{aligned}
    |\widehat{x}(t_q^-)-\xi - x(t_q^-)| & \leq \beta(|\widehat{x}(t_p^-)-x(t_p^-)|,t_q-t_p) \\
    & < |\widehat{x}(t_p^-)-x(t_p^-)|, 
    \end{aligned}
 \]
 i.e.\ the correction term $\xi$ indeed reduces the estimation error at the end of the interval in comparison to the estimation error at the beginning of the interval.
\end{Theorem}
\begin{Proof}
Let $\widehat{x}^\new$ be the (virtual) solution of the system copy \eqref{eq:system_copy_on_[tp,tq)} with corrected initial value $\widehat{x}^\new(t_p^-) = \widehat{x}(t_p^-) - \xi^\myleft$. 
   By construction $\widehat{x}^\new(t_p^-)-x(t_p^-) = e(t_p^-)-\xi^\myleft\in\cN_p^q$ hence $y=\widehat{y}^\new$ on $[t_p,t_q)$ and therefore, for all $t\in[t_p,t_q)$,
    \[
    |\widehat{x}^\new(t^+)-x(t^+)| \leq \beta(|\widehat{x}^\new(t_p^-)-x(t_p^-)|,t-t_p).
 \]
 Note that $e=\widehat{x}-x$ as well as $e^\new:=\widehat{x}^\new - x$ are solutions of the homogenous error DAE \eqref{eq:errorDAE}, in particular
 \[
    e(t_q^-) - e^\new(t_q^-) = \Phi_p^q (e(t_p^-)-e^\new(t_p^-)) = \Phi_p^q \xi^\myleft = \xi,
 \]
 or, in other words,
 \[
    \widehat{x}(t_q^-) - \xi = \widehat{x}^\new(t_q^-)
 \]
Finally, by construction $\xi^\myleft\in{\cN_p^q}^\bot$ and therefore, by Pythagoras' Theorem,
 \[\begin{aligned}
     |\widehat{x}(t_p^-)-x(t_p^-)|^2 & = |\underbrace{\widehat{x}^\new(t_p^-) - x(t_p^-)}_{\in\cN_p^q} +\xi^\myleft|^2 \\
     & = |\widehat{x}^\new(t_p^-) - x(t_p^-)|^2 + |\xi^\myleft|^2\\
     & \geq |\widehat{x}^\new(t_p^-)-x(t_p^-)|^2.
   \end{aligned}
 \]
Altogether we have:
\[
\begin{aligned}
  |\widehat{x}(t_q^-)-\xi - x(t_q^-)| &= |\widehat{x}^\new(t_q^-) - x(t_q^-)| \\
  & \leq \beta(|\widehat{x}^\new(t_p^-) - x(t_p^-)|,t_q-t_p)\\
   &\leq \beta(|\widehat{x}(t_p^-) - x(t_p^-)|,t_q-t_p)
  \end{aligned}
\]
which is the desired estimate.
\end{Proof}
\begin{Remark}\label{rem:correction_xi^left}
   The proof of Theorem \ref{thm:xi_correction} reveals that by applying the correction $\xi^\myleft$ at the beginning of the interval the output of the system copy is then identical to the output of the original system. However, for the observer design it is \emph{not} necessary to rerun the system copy (in particular storing the whole input signal over the interval $[t_p,t_q)$), because we just propagate $\xi^\myleft$ via the homogenous error dynamics \eqref{eq:errorDAE} which is independent of the in- and output. This actually allows us to calculate the error correction for an arbitrary future; this fact can be utilized to deal with time delays due to computation times, see the discussion in Section~\ref{sec:implementation}.
\end{Remark}

\section{Estimation errors and asymptotic convergence}\label{sec:convergence}

In theory, it is possible to determine the observable part exactly from the output, however, in practice one can only get approximations. Nevertheless, these approximations may be as accurate as desired (e.g.\ by choosing appropriate gains in a Luenberger observer). Similar as in \citep{TanwTren17a} we therefore make the following assumption about the ability to approximate the observable part to any given accuracy:
\begin{Assumption}\label{ass:hat_zk}
   For each mode $k$ of the switched DAE \eqref{eq:swDAE} and a given $\eps_k>0$, there exists an estimator $\widehat{z}_k=\widehat{\fO}_k(y^e_{[t_k,t_{k+1})})$ such that
   \begin{equation}\label{eq:epsReq}
      |\widehat{z}_k - z_k| \leq \eps_k |z_k|,
   \end{equation}
   where $z_k=\fO_k(y^e_{[t_k,t_{k+1})})$ is the ideal estimator of the observable part on $[t_k,t_{k+1})$ as given in Section~\ref{sec:zk}.
\end{Assumption}

Under Assumption \ref{ass:hat_zk}, the state estimation correction in \eqref{eq:system_copies_with_xi_i} for the interval $[t_{p_i},t_{q_i})$ is given by
\begin{equation}\label{eq:xi_definition}
   \xi_i:=\Phi_{p_i}^{q_i}\cO_{p_i}^{q_i-1}\widehat{\mathbf{z}}_{p_i}^{q_i-1},
\end{equation}
where $\widehat{\mathbf{z}}_{p_i}^{q_i-1} = (\widehat{z}_{p_i}/\widehat{z}_{p_i+1}/\cdots/z_{q_i-1})$.

As detailed in~\ref{sec:localObsComp}, the observable component $z_k$ of the estimation error $e=\widehat{x}-x$ on the interval $[t_k,t_{k+1})$ is composed of the two components $z_k^\diff$ and $z_k^\imp$, where the former is obtained from the continuous output mismatch $y^e$ on $(t_k,t_{k+1})$ and the latter is obtained from the impulsive mismatch $y^e[t_k]$. The estimation of $z_k^\diff$ can be reduced to the classical state estimation problem for non-switched linear ODEs and there are many methods to do that. The only non-standard aspect here is that we have to obtain the state-estimation at the beginning of the interval $(t_k,t_{k+1})$ and not (as usual) at the end of the interval.
This does not pose any serious problems, as we can use a standard Luenberger observer on the interval $(t_k,t_{k+1})$ to get an estimate at the end of the interval and then propagate this estimate back in time. Since the (homogeneous) ODE dynamics are known as well as the length of the interval, we can ensure the desired estimation accuracy at the beginning of the interval by increasing the accuracy of the estimate at the end of the interval.\footnote{In fact, 
consider an observable LTI system $\dot z = A z$, $y= Cz$ over the interval $[0,T]$, with the estimator $\hat z = (A-LC)\hat z + Ly$, $\hat z(0) = 0$. For every $\epsilon> 0$, there exists $L$ such that $\vert \hat z(T) - z(T) \vert \le \eps \vert z(0)\vert$. Choose $\hat z_0^* = e^{-AT}\hat z(T)$, and $\epsilon \le \delta/\|e^{-AT}\|$ for some desired $\delta > 0$, then it is easily verified that $\vert z(0) - z_0^* \vert \le \delta \vert z(0) \vert$.
}
Another (more sophisticated) way of obtaining such estimates is by the use of ``back-and-forth observer'' as presented in \citep{ShimTanw12}, however, this requires the storage of the output over the whole interval $(t_k,t_{k+1})$.

The estimation accuracy for $z_k^\imp$ is actually concerned with the measurement accuracy of the impulsive part $y^e[t_k]$, i.e.\ on how well Dirac impulses and their derivatives can be measured in practice, see \citep{TanwTren17a} for details.

Assumption~\ref{ass:hat_zk}, together with Assumption~\ref{ass:uniform_localdetec}, provide all the ingredients we need for obtaining converging state estimates.


\begin{Theorem}\label{thm:convergence}
   Consider the switched DAE \eqref{eq:swDAE} satisfying the uniform local detectability Assumption~\ref{ass:uniform_localdetec}, and the local estimation accuracy Assumption~\ref{ass:hat_zk}.
   For the $\alpha$ given in \eqref{eq:betai_alpha}, choose $\eps_k$, $k\in\N$, (depending on $\alpha$) such that
    \begin{equation}\label{eq:epsChoice}
    \boxed{c_i \eps^{\max}_i \leq \widehat{\alpha}-\alpha}
 \end{equation}
 for some $\widehat \alpha \in (\alpha,1)$,  where
   \[
      c_i := \|\Phi_{p_i}^{q_i}\cO_{p_i}^{{q_i}-1}\| \left\|\begin{bmatrix} Z_{p_i}^\top \\ Z_{p_i+1}^\top \Phi_{p_i}^{p_i+1} \\ \vdots \\ Z_{q_i-1}^\top \Phi_{p_i}^{q_i-1} \end{bmatrix}\right\|
   \]
   and 
   \[
     \eps_i^{\max}:=\max\setdef{\eps_k}{p_i\leq k \leq q_i-1}. 
   \]
Then the observer given by the system copies \eqref{eq:system_copies_with_xi_i}, with error corrections $\xi_i$ in \eqref{eq:xi_definition} and the estimate $\widehat z_k$ chosen to satisfy \eqref{eq:epsReq} for $\eps_k$ specified in \eqref{eq:epsChoice}, results in
\[
    \widehat{x}(t^+)\to x(t^+)\text{ as }t\to\infty,
\]
i.e.\ the observer achieves asymptotic estimation of the state.
\end{Theorem}
\begin{Proof}
   From Theorem~\ref{thm:xi_correction} we know that for each detectable interval $[t_{p_i},t_{q_i})$ the ideal correction $\xi_i^\ideal=\Phi_p^q\cO_p^{q-1}\mathbf{z}_p^{q-1}$ ensures
   \[
   \begin{aligned}
      |\widehat{x}_i(t_{q_i}^-) - \xi_i^\ideal - x(t_{q_i}^-)| & \leq \beta_i(|\widehat{x}_i(t_{p_i}^-) - x(t_{p_i})|,t_{q_i}-t_{p_i}) \\
      & \leq \alpha |\widehat{x}_i(t_{p_i}^-) - x(t_{p_i})|.
    \end{aligned}
   \]
   Hence for the actual correction term $\xi_i = \Phi_p^q\cO_p^{q-1}\widehat{\mathbf{z}}_p^{q-1}$, we have
   \[
   \begin{aligned}
      |\widehat{x}_{i+1}(t_{q_i}^-) - x(t_{q_i}^-)| \leq  \alpha |\widehat{x}_i(t_{p_i}^-) - x(t_{p_i})| + |\xi_i^\ideal - \xi_i|.
      \end{aligned}
   \]
   By assumption, $|\widehat{z}_k-z_k|\leq \eps_k |z_k|$ and since $z_k = Z_k^\top e(t_k^-) = Z_k^\top \Phi_{p_i}^k e(t_{p_i}^-)$ for any $k$ with $p_i\leq k\leq q_i-1$, we have
   \[\begin{aligned}
      & |\xi_i^\ideal - \xi_i| = |\Phi_{p_i}^{q_i}\cO_{p_i}^{{q_i}-1}(\widehat{\mathbf{z}}_{p_i}^{{q_i}-1} - \mathbf{z}_{p_i}^{{q_i}-1})|\\
      &\quad \leq \|\Phi_{p_i}^{q_i} \cO_{p_i}^{{q_i}-1}\| \eps_i^{\max} \left\|\begin{bmatrix} Z_{p_i}^\top \\ Z_{p_i+1}^\top \Phi_{p_i}^{p_i+1} \\ \vdots \\ Z_{q_i-1}^\top \Phi_{p_i}^{q_i-1} \end{bmatrix}\right\|  |e(t_{p_i}^-)|\\
      & \quad =c_i \eps^{\max}_i |\widehat{x}_i(t_{p_i}^-)-x(t_{p_i}^-)|.
   \end{aligned}
   \]
 Altogether we have
   \[
      |\widehat{x}_{i+1}(t_{q_i}^-) - x(t_{q_i}^-)| \leq \widehat{\alpha} |\widehat{x}_{i}(t_{p_i}^-) - x(t_{p_i}^-)|,
   \]
   i.e.\ on each detectability interval $[t_{p_i},t_{q_i}$) the estimation error $\widehat{x}-x$ decreases uniformly by a factor $\widehat{\alpha}$ and the same proof technique as in Proposition~\ref{prop:uniform_local_implies_global} shows asymptotic convergence.
\end{Proof}

\section{Observer implementation}\label{sec:implementation}
We have implemented the observer in Matlab and show in the next section the simulation results for the academic Examples \ref{ex:swODEdim3} and \ref{ex:impulses_important}. Before presenting the simulation results, we would like to discuss some implementation issues, in particular, which calculations can be carried out offline and how to treat the necessary computation times. The overall structure of the observer is given in Algorithm~\ref{alg:observer}.

\begin{algorithm}[!hbt]\caption{Observer for detectable switched DAEs}\label{alg:observer}
\small
\KwData{Modes $(E_k,A_k,B_k,C_k,D_k)$, $k=0,1,2,\ldots$\\
  switching times~$t_k$, $k=1,2,3,\ldots$, $t_0:=0$\\
  update-time indicies $q_i$, $i=0,1,2,\ldots$, $q_{-1}:=0$\\
  access to input $u$ and output $y$}
\KwResult{State-estimation $\widehat{x}$}
\textsc{Initialization (offline):}\\
\ForAll{modes $k$}{
   Calculate $\Pi_k$, $A^\diff_k$, $E^\imp_k$, $C^\diff_k$ as in Def.~\ref{def:proj}\\
   Calculate $\me^{A^\diff \tau_k}\Pi_k$ and $\me^{-A^\diff_k \tau_k}$\\
   Calculate $W_k$, $Z_k$, $Z_k^\diff$, $Z_k^\imp$, $U_k^\obs$ as in App.~\ref{sec:localObsComp}\\
   Calculate $S_k^\diff$, $R_k^\diff$, $U^\imp_k$ as in \eqref{eq:zsys}, \eqref{eq:Ukimp}\\
   Choose $L_k$ s.\ t.\ $S_k^\diff-R_k^\diff L_k$ is ``sufficiently'' Hurwitz (in view of Assumption~\ref{ass:hat_zk})\\
}
\textsc{Run observer on detectable intervals $[p_i,q_i)$:}\\
\ForAll{$i\in\N$}{
$p:=q_{i-1}$,\quad $q:=q_i$\\
\textsc{Get local estimation data (online):}\\
\ForAll{$k = p, p +1 , \ldots, q-1$}{
   Run system copy \eqref{eq:system_copies_with_xi_i} with input $u$\\
   \textbf{Run Luenburger observer} with gain $L_k$ for \eqref{eq:zsys} on $(t_k,t_{k+1})$ with output-injection $y^e\!=\!\widehat{y}\!-\!y$\\
   $\to$ estimation of $\mathbf{z}^\diff_k$ on $(t_k,t_{k+1})$\\
   \textbf{Estimate impulse differences} $\widehat{y}[t_k]-y[t_k]$\\
   $\to$ estimation of $\boldsymbol{\eta}_k$\\
   Calculate $\widehat{z}^\diff_k$, $\widehat{z}^\imp_k$ according to \eqref{eq:hatzkdiff},\eqref{eq:zkimp}\\
   Calculate $\widehat{z}_k$ via \eqref{eq:calculate_zhat}\\
}
\textsc{\mbox{Combine local information backwards (offline)}}\\
   $\widehat{\mu}_{q-1} := \widehat{z}_{q-1}$\\
   \ForAll{$k = q-2,q-3,\ldots,p$ (backwards)}{
     Calculate subspaces $\cN^{q}_{k}$ recursively via \eqref{eq:DAESeqUnobs}\\
     Choose $\Theta^{q}_k$, $M^{q}_k$, $U^q_k$ via \eqref{eq:Thetaqk_def}, \eqref{eq:Mqk_def}, \eqref{eq:Uqk_def}\\
     Calculate $\widehat{\mu}_{k}$ from $\widehat{z}_k$, $\widehat{\mu}_{k+1}$ via \eqref{eq:muk_recursion}\\
  }
  $\xi^\myleft := M^q_p \widehat{\mu}_p$\\
  \textsc{Propagate correction forward (offline)}\\
  $\xi_i \leftarrow \xi^\myleft$\\
  \ForAll{$k = p+1, p+2, \ldots, q-1$ (forward)}{
     Calculate $\xi_i \leftarrow \me^{A^\diff_k \tau_k} \Pi_k \xi_i$
  }
  Update state-estimation $\widehat{x}(t_{q_i}^-)\leftarrow \widehat{x}(t_{q_i}^-) - \xi_i$\\
}
\end{algorithm}

In the initialization phase, certain matrices and subspaces are calculated for each individual subsystem; in particular, a decomposition into unobservable and observable states is carried out. It should be noted that, in practice and in our setup, these calculations have to be carried out only for finitely many modes. In fact, it suffices to carry out the calculation for all modes occurring in the next detectability interval (and these calculations can be done in parallel to running the system copy and collecting the corresponding measurements). For a suitable choice of the Luenberger gain it is necessary to know (at least some bounds of) the values $c_i$ in \eqref{eq:epsChoice} and consequently (bounds on) the necessary estimation accuracies $\eps_k$ as in Assumption~\ref{ass:hat_zk}. In case of a periodic switching signal one may adapt the Luenberger gain in each iteration until the estimation error is sufficiently small such that convergence of the state-estimator is guaranteed.

In the online (or synchronous) phase of the observer, a system copy must be simulated (driven by the current input) and its output must be compared with the actual output of the system. Although the output of the original system may contain Dirac impulses at arbitrary times, for the observer only the Dirac impulses \emph{at the switching times} are compared with the predicted Dirac impulses of the system copy at the switching times. All other Dirac impulses are induced by discontinuities in the input (and are independent of the current state), hence (at least in theory) they are identical for the system copy and the original system and do not provide any additional information for the state estimation problem. Hence the impulse measurement needs only be active around the switching times. For estimating the observable part between the switching times, one could either run a classical observer (for the system \eqref{eq:zsys} with the desired output $y^e=\widehat{y} - y$) synchronously to the system without the need to store the measures output difference $y^e$; however, as the dimension of the output is usually low, it may also be feasible to store the whole trajectory $y^e$ and carry out some more sophisticated estimation procedure offline.

Finally, after the local observability data is obtained, it must be combined in a suitable way to obtain the impulsive update $\xi$ for the state estimation $\widehat{x}$. Although all the involved matrices can be computed offline, the actual calculations can only be carried out after the last estimate $\widehat{z}_{q_i-1}$ is obtained, hence some unavoidable processing time $\Delta>0$ is required to compute $\xi$. However, the effect of the processing time can be entirely compensated as follows: For a generic detectability interval $[t_p,t_q)$, assume that an upper bound $\Delta> 0$ is known for the time required to calculate $\xi$. Furthermore, we assume that $[t_p,t_q-\Delta)$ is still a detectable interval in the sense of Definition~\ref{def:localDet} (this is always the case for sufficiently small $\Delta$). In particular, $\xi^\myleft$ is an arbitrarily good estimate of the projection of $e(t_p^-)$ on the unobservable space $\cN^q_p$. Now, we just propagate forward $\xi^\myleft$ with the matrix $\Phi^q_p$ to get a good estimate of $e(t_q^-)$, and we can update $\hat x$ at the correct time. The key observation is that once we have obtained $\xi^\myleft$, we can freely chose the update time (i.e.\ how far we propagate forward the error correction) without loosing any accuracy.

\section{Simulations}\label{sec:simulations}
\subsection{Simulation of Example~\ref{ex:swODEdim3}} \label{ex:simSwODEdim3}

Consider the switched ODE given in Example~\ref{ex:swODEdim3} with the periodic switching signal. It was already shown in Example~\ref{ex:detODEdim3} that, due to periodicity assumption on $\sigma$, this system is uniformly interval-detectable. To implement the proposed observer, we run the system copy \eqref{eq:system_copies_with_xi_i} on the intervals $[3i,3i+3)$, $i \in \N$, and apply the correction term at $t_{q_i} = 3i+3$. The correction terms are obtained by
\[
\xi_{i} = \me^{A_1\tau_{3i+1}} M_{3i}^{3i+3} \widehat \mu_{3i}, \quad i \in \N,
\]
where $\tau_{3i+1} = 1$, for each $i$, and $\widehat \mu_{3i}$ is computed from the estimates of the observable states of individual subsystems: $\widehat z_{3i}$ and $\widehat z_{3i+2}$. We recall that $z_{3i+1}$ is an empty vector because the output over the interval $[3i+1,3i+2)$ is zero for each $i$ and nothing can be deduced about the state. The values of $\widehat z_{3i}$ and $\widehat z_{3i+2}$ are obtained by running a Luenberger observer (with gain $L=1$) for the $x_1$-dynamics over the intervals $[3i,3i+1)$ and $[3i+2,3i+3)$, respectively. 
We use the later to first compute
\[
\widehat\mu_{3i+1} = \cF_{3i+1}^{3i+3} \begin{pmatrix} \widehat z_{3i+1} \\ \widehat z_{3i+2} \end{pmatrix} = \Theta_{3i+1}^\top \me^{-A_{3i+1} \tau_{3i+1}} Z_{3i+2} \widehat z_{3i+2}
\]
where $\Theta_{3i+1}^\top = \left(-1/\sqrt{2}, -1/\sqrt{2}, 0\right)$, and $Z_2 = \left(1/0/0\right)$. This leads to 
\[
\widehat \mu_{3i} = \cF_{3i}^{3i+3} \begin{pmatrix} \widehat z_{3i} \\ \widehat \mu_{3i+1} \end{pmatrix} = \begin {pmatrix} \widehat z_{3i} \\ \Theta_{3i}^\top M_{3i+1} \widehat \mu_{3i+1} \end{pmatrix}
\]
where $\Theta_{3i+1}^\top = \left(1, 1, 0\right)$ and $M_{3i+1}^{3i+3^\top} = \begin{smallbmatrix} 1 & 0 & 0\\ 0 & 1 & 0\end{smallbmatrix}$.

\begin{figure}[t]
  \centering
  \includegraphics[width=0.485\textwidth]{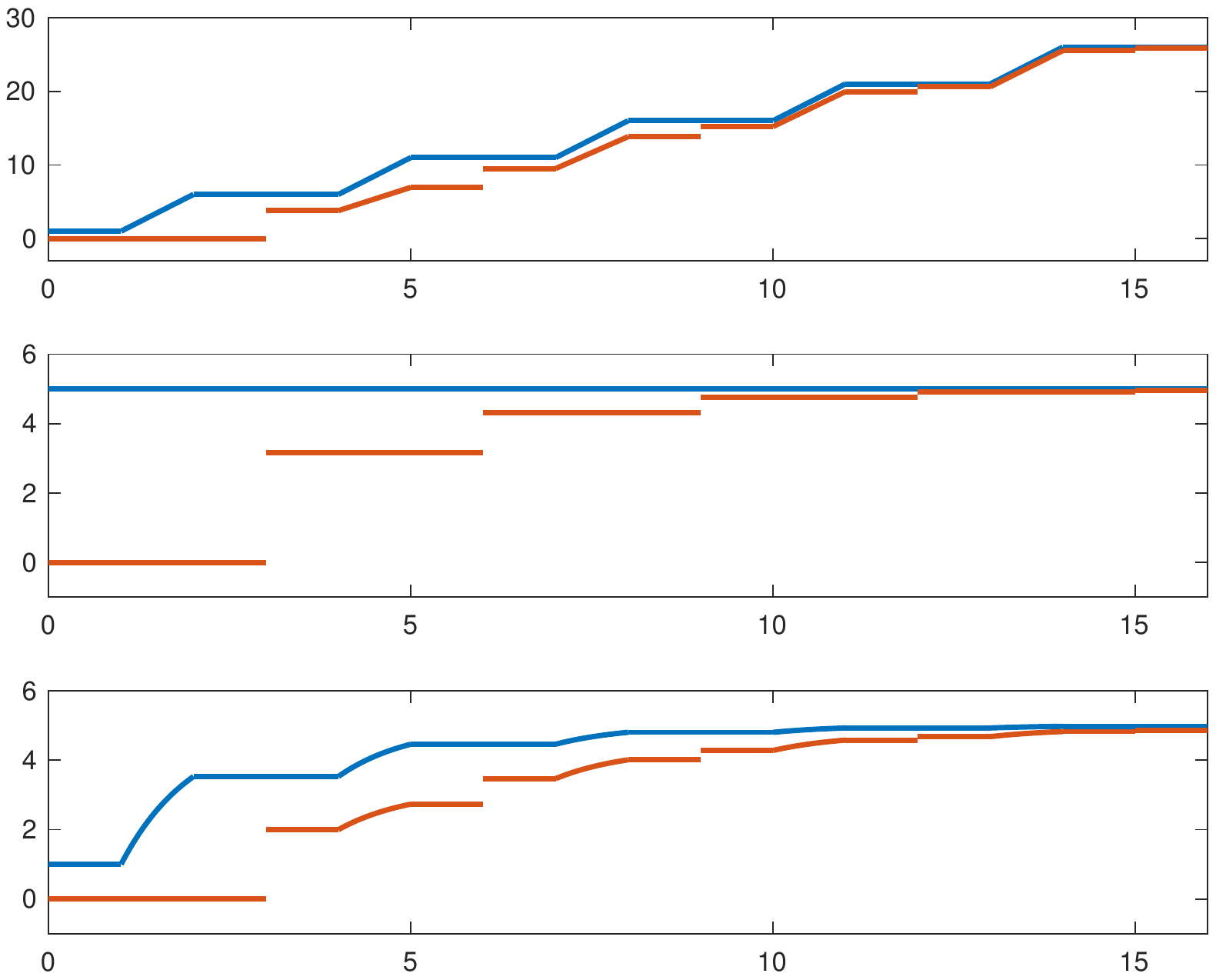}
  \caption{State estimation, $x_1$ (blue) and $\widehat{x}_1$ (red) top figure; $x_2$ (blue) and $\widehat{x}_2$ (red) middle figure; $x_3$ (red) and $\widehat{x}_3$ (red) bottom figure.}\label{fig:errors}
\end{figure}

The results of the simulation are reported in Figure~\ref{fig:errors}. 

It is observed that whenever a correction is applied at $t_{q_i} = 3i+3$, the estimation error decreases and the convergence to zero is achieved asymptotically.

\subsection{Simulation of Example~\ref{ex:impulses_important}}

We now implement our observer on the system given in Example~\ref{ex:impulses_important} where one of the subsystem is a DAE. As already discussed above the presence and evaluation of the occurring Dirac impulses in the output are crucial for the state estimation.

The system is detectable on the intervals $[2i,2i+2)$, for $i \in \N$ and so we run the system copy \eqref{eq:system_copies_with_xi_i} on the intervals $[2i,2i+2)$, and apply the correction term at $t_{q_i} = 2i+2$. The correction term, for each $i \in \N$, is obtained by
\[
\xi_i = \me^{A_{2i+1}^\diff \tau_{2i+1}} \Pi_{2i+1}\me^{A_{2i}^\diff \tau_{2i}} \Pi_{2i} M_{2i}^{2i+2} \widehat \mu_{2i},
\]
where $\tau_{2i}= \tau_{2i+1} =1$ and
\begin{gather*}
A_{2i}^\diff = A_{2i} = \begin{smallbmatrix} 0 & 0 & 1 & 0\\0 & 0 & 0 & 0\\0 & 0 & 0 & 0 \\ 0 & 0 & 1 & -1\end{smallbmatrix}, \quad 
A_{2i+1}^\diff = \begin{smallbmatrix} 0 & 0 & 0 & 0\\0 & 0 & 0 & 0\\0 & 0 & 0 & 0 \\ 0 & 0 & 1 & -1\end{smallbmatrix},\\
\Pi_{2i} = I_{4 \times 4}, \quad \Pi_{2i+1} = \begin{smallbmatrix} 0 & 0 & 0 & 0 \\ 0 & 0 & 0 & 0 \\ 0 & 0 & 1 & 0 \\ 0 & 0 & 0 & 1\end{smallbmatrix}.
\end{gather*}

\begin{figure}[t]
  \includegraphics[width=0.485\textwidth]{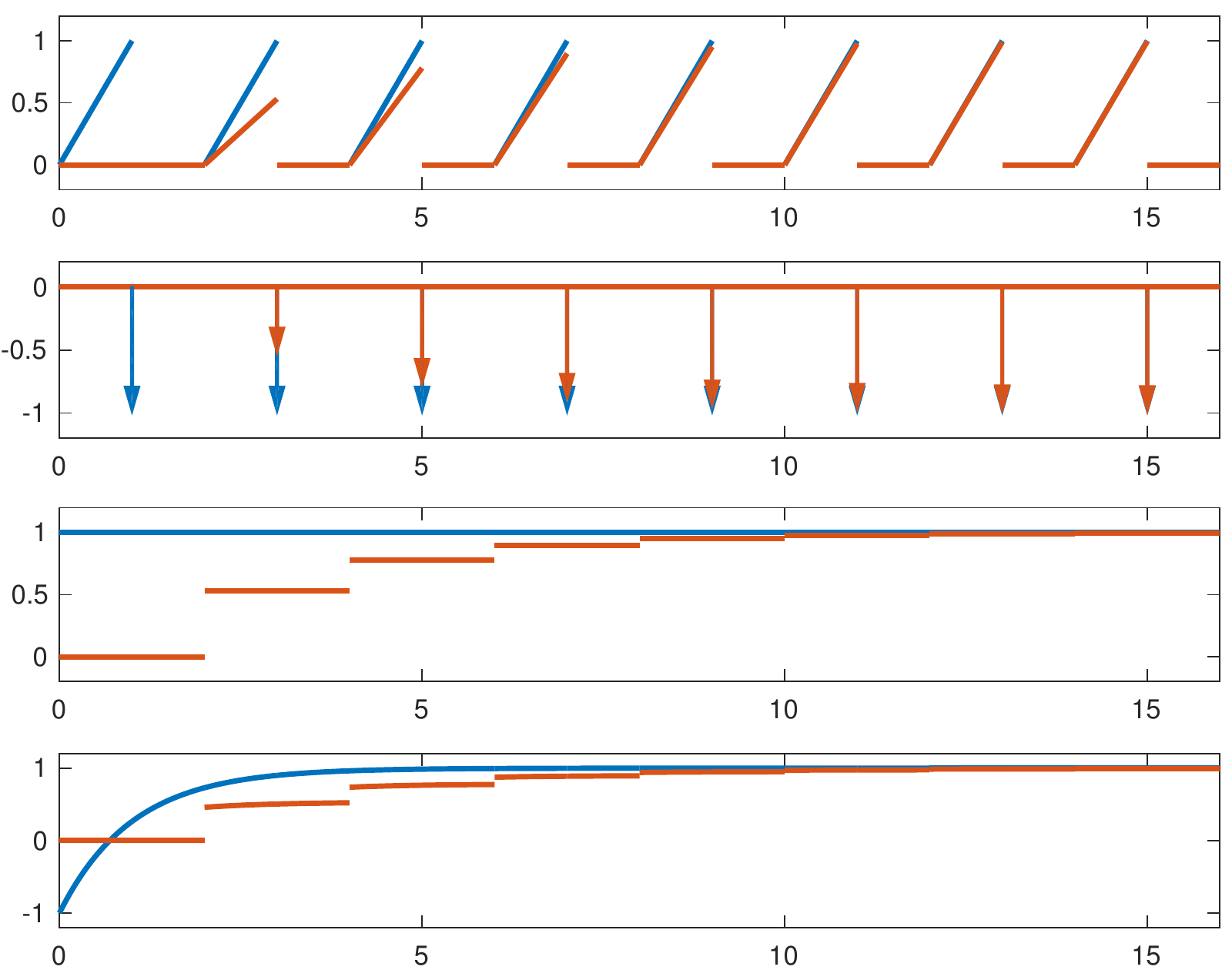}
  \caption{State estimation, $x_1$ (blue) and $\widehat{x}_1$ (red) in top figure; $x_2$ (blue) and $\widehat{x}_2$ (red) with Dirac impulses (shown as arrows) in second figure from top (the Dirac impulses are also visible in the output); $x_3$ (blue) and $\widehat{x}_3$ (red) in third figure; $x_4$ (blue) and $\widehat{x}_4$ (red) in bottom figure.}\label{fig:errorEx2}
\end{figure}
 
The term $\widehat \mu_{2i} \in \R$ is computed from the estimates of the observable states of individual subsystems: $\widehat z_{2i}$ and $\widehat z_{2i+1}$. Because $y = 0$ on $[2i,2i+1)$, and the corresponding subsystem is an ODE, we set $\widehat z_{2i}$ to be the empty vector (nothing can be concluded from the output). Also, due to the structure of the second subsystem, the only observable information is due to impulses in the output, so that $\widehat z_{2i+1} = \widehat z_{2i+1}^\imp \in \R$ because
\[
O_{2i+1}^\diff = 0_{4 \times 4}, \quad O_{2i+1}^\imp = \begin{smallbmatrix} 1 & 0 & 0 & 0 \\ 0 & 0 & 0 & 0 \\ 0 & 0 & 0 & 0 \\ 0 & 0 & 0 & 0 \end{smallbmatrix}.
\]
To compute $\widehat \mu_{2i}$, we have $\widehat \mu_{2i+1} = \widehat z_{2i+1}$, and 
\[
\widehat \mu_{2i} = \cF_{2i}^{2i+2} \widehat \mu_{2i+1} = \Theta_{2i}^\top \me^{-A_{2i}^\diff \tau_{2i}} Z_{2i+1} \widehat z_{2i+1}
\]
where $\Theta_{2i}^\top = \left(\frac{-1}{\sqrt{2}}, 0 , \frac{-1}{\sqrt{2}} , 0\right) $, and $Z_{2i+1}^\top = (1,0,0,0)$. Finally, we compute $M_{2i}^{2i+2^\top} = \begin{smallbmatrix} \frac{1}{\sqrt{2}}, 0 , \frac{1}{\sqrt{2}} , 0 \end{smallbmatrix} $ and use it along with $\widehat \mu_{2i}$ to compute the correction terms $\xi_{2i+2}$.

Since in the simulation both the system and the system copy are simulated, we would be able to obtain $\widehat{z}^\imp_{2i+1}$ without estimation error (i.e.\ $\eps_{2i+1}=0$ in Assumption~\ref{ass:hat_zk}) and then already after one correction we would have a perfect state-estimation. Therefore, we introduced some artificial random noise when ``measuring'' $y[t_{2i+1}]$, so that \eqref{eq:epsReq} is only satisfied with $\eps_{2i+1}=0.1$.

The results of the simulation are reported in Figure~\ref{fig:errorEx2} where we see that the estimation error converges to zero in each of the state components. In particular, we see corrections in the magnitude of the impulses in state component $x_2$.

\section{Conclusion}
We have studied the notion of detectability for switched DAEs which allows us to consider the problem of observer design under relaxed assumptions on system dynamics compared to the existing works. A novel estimation algorithm is proposed which relies on propagating backward and forward the correction terms obtained by processing the measured outputs and inputs. Rigorous convergence analysis of estimation error for the proposed algorithm is carried out and the results are illustrated by studying two academic examples with simulations.

\section*{Appendix}
\appendix

\section{Properties of a matrix pair $(E,A)$} \label{sec:basicProp}
A very useful characterization of regularity is the following well-known result.
\begin{Proposition}[Regularity and quasi-Weierstra\ss\  form]
A matrix pair $(E,A)\in\R^{n\times n}\times\R^{n\times n}$ is regular if, and only if, there exist invertible matrices $S,T\in\R^{n\times n}$ such that 
\begin{equation}\label{eq:QWF}
(SET,SAT) = \left( \begin{bmatrix} I & 0  \\  0  &  N \end{bmatrix}, \begin{bmatrix} J  & 0 \\  0  & I \end{bmatrix} \right) ,
\end{equation}
where $J \in \R^{n_1\times n_1}$, $0\leq n_1\leq n$,  is some matrix and $N\in\R^{n_2 \times n_2} $, $n_2:=n-n_1$, is a nilpotent matrix.
\end{Proposition}

We call \eqref{eq:QWF} a \emph{quasi-Weierstrass form} of $(E,A)$ following \cite{BergIlch12a}; therein it also shown how to easily obtain \eqref{eq:QWF} via the Wong-sequences \citep{Wong74}.

\begin{Definition}\label{def:proj}
 Consider the regular matrix pair $(E,A)$ with corresponding quasi-Weierstra\ss\ form \eqref{eq:QWF}.
 The \emph{consistency projector} of $(E,A)$ is given by
\[
   \Pi = T \begin{bmatrix} I & 0\\ 0 & 0 \end{bmatrix}T^{-1}.
\]
Furthermore, let
\[\begin{aligned}
    A^\diff &:= T \begin{bmatrix} J & 0 \\ 0 & 0 \end{bmatrix} T^{-1},\\
    E^\imp &:= T \begin{bmatrix} 0 & 0 \\ 0 & N \end{bmatrix} T^{-1}.
    \end{aligned}
\]
Finally, if also an output matrix $C$ is considered let
\[
    C^\diff := C \Pi_{(E,A)}.
\]
\end{Definition}

To see the utility of the matrices introduced in Definition~\ref{def:proj}, consider the problem of finding a trajectory $x$ which solves the initial-trajectory problem (ITP)
\begin{subequations}\label{eq:ITP}
  \begin{align}
    x_{(-\infty,0)}&=x^0_{(-\infty,0)}\label{eq:ITPa}\\
    (E\dot{x})_{[0,\infty)} &= (Ax)_{[0,\infty)}, \label{eq:ITPb}
  \end{align}
\end{subequations}
in some appropriate sense.

\begin{Lemma}[Role of consistency projector, {\cite[Thm.~4.2.8]{Tren09d}}]\label{lem:consProj}
   Consider the ITP \eqref{eq:ITP} with regular matrix pair $(E,A)$ and with arbitrary initial trajectory $x^0\in(\Dpwsm)^n$. There exists a unique solution $x\in(\Dpwsm)^n$ and
   \[
      x(0^+)=\Pi_{(E,A)} x(0^-).
   \]
\end{Lemma}

\begin{Lemma}[{\cite[Cor.~5]{TanwTren10}}]\label{lem:impulses}
Consider the ITP \eqref{eq:ITP} with regular matrix pair $(E,A)$ and the corresponding $E^\imp$ matrix.
For the unique solution $x\in(\Dpwsm)^n$, it holds that
  \begin{equation}\label{eq:impSol}
     x[0] = -\sum_{j=0}^{n-2}(E^\imp)^{j+1} x(0^-) \delta_0^{(j)},
  \end{equation}
  where $\delta_0^{(j)}$ denotes the $j$-th (distributional) derivative of the Dirac-impulse $\delta_0$ at $t=0$.
\end{Lemma}

\begin{Lemma}\label{lem:diffProj}
  For any regular matrix pair $(E,A)$ and output matrix $C$, the following implication holds for all continuously differentiable $(x,y)$:
  \[\left.\begin{aligned}
     E\dot{x} &= Ax,\\
     y &= Cx
     \end{aligned}\right\}\ \Rightarrow\ 
     \left\{\begin{aligned}
        \dot{x} &= A^\diff x,\\
        y &= C^\diff x
        \end{aligned}\right.
  \]
  In particular, any classical solution $x$ of $E\dot{x}=Ax$ satisfies
  \[
     x(t) = \me^{A^\diff t} x(0),\quad t\in\R. 
  \]
\end{Lemma}

\section{Output-to-State Mappings}

\subsection{Observable component of a subsystem}\label{sec:localObsComp}

The local unobservable space \eqref{eq:Wk} is given by
\[
   \cW_k = \Pi_k^{-1}\ker O^\diff_k\ \cap\ \ker O^\imp_k,
\]
where
\begin{equation}\label{eq:Odiff_Oimp}
\begin{aligned}
O_k^\diff & :=[C_k^\diff / C_k^\diff A_k^\diff / \cdots / C_k^\diff (A_k^\diff)^{n-1}], \\
O_k^\imp &:=[C_k E_k^\imp / C_k (E_k^\imp)^2 / \cdots / C_k (E_k^\imp)^{n-1}].
\end{aligned}
\end{equation}
In other words, $\ker O_k^\diff$ denotes the unobservable space of the ODE $\dot{e}=A_k^\diff e$, $y^e=C_k^\diff e$, 
and $\ker O_k^{\imp}$ denotes the impulse unobservable space in the sense that $y^e[t_k] = 0$ implies $e(t_k^-) \in \ker O_k^{\imp}$.

We may now write
\[
e(t_k^-) = W_kw_k + Z_k z_k,
\]
where $\im W_k=\cW_k$ and $\im Z_k = \cW_k^\bot$ and $W_k,Z_k$ are orthonormal matrices. Here $z_k$ determines the projection of $e(t_k^-)$ onto the subspace $\cW_k^\bot$. The latter can further be decomposed as
\[
   \cW_k^\bot = \im (O_k^\diff \Pi_k)^\top + \im {O_k^\imp}^\top
\]
Let $Z_k^\diff $, and $Z_k^\imp $ be the orthonormal matrices such that
\[\begin{aligned}
   & \im Z_k^\diff = \im\left({O^\diff_k}^\top\right), & \quad & z_k^\diff := {Z_k^{\diff}}^\top \Pi_k e(t_k^-) \\
   & & \quad & \phantom{z_k^{\diff}:} = {Z_k^{\diff}}^\top e(t_k^+), \\
   & \im Z_k^\imp = \im\left({O^\imp_k}^\top\right), & \quad & z_k^\imp := {Z_k^{\imp}}^\top e(t_k^-).
\end{aligned}
\]
The motivation for introducing the components $z_k^\diff$ and $z_k^\imp$ is that they can be estimated using the output measurements on the interval $[t_k,t_{k+1})$.
To express the vector $z_k$ in terms of these components, we introduce the matrix $U_k^\obs$ such that
\[
Z_k = \begin{bmatrix} \Pi_k^\top Z_k^\diff & Z_k^\imp \end{bmatrix} U_k^\obs.
\]
Such a matrix $U_k^\obs$ always exists because
\[
\begin{aligned}
\im Z_k = \cW_k^\bot & = (\Pi_k^{-1}(\ker O_k^\diff))^\bot + (\ker O_k^\imp)^\bot \\
& = \Pi_k^\top \im Z_k^\diff + \im Z_k^\imp\\
& = \im \begin{bmatrix} \Pi_k^\top Z_k^\diff & Z_k^\imp \end{bmatrix}.
\end{aligned}
\]
It then follows that
\[
\begin{aligned}
z_k = Z_k^\top e(t_k^-) & =  {U_k^{\obs}}^\top \begin{bmatrix} Z_k^{\diff^\top} \Pi_k \\ Z_k^{\imp^\top} \end{bmatrix} e(t_k^-) \\
& ={U_k^{\obs}}^\top \begin{bmatrix} Z_k^{\diff^\top} e(t_k^+) \\ Z_k^{\imp^\top} e(t_k^-) \end{bmatrix} = {U_k^{\obs}}^\top \begin{bmatrix} z_k^\diff \\ z_k^\imp \end{bmatrix}.
\end{aligned}
\]
If only estimates $\widehat{z}^\diff_k$ and $\widehat{z}^\imp_k$ of $z^\diff_k$ and $z^\imp_k$ are available, we therefore obtain an estimate of $z_k$ as follows:
\begin{equation}\label{eq:calculate_zhat}
   \widehat{z}_k = {U_k^{\obs}}^\top \begin{bmatrix} \widehat{z}_k^\diff \\ \widehat{z}_k^\imp \end{bmatrix}.
\end{equation} 
Next, we specify how to write $z_k^\diff$ and $z_k^\imp$ in terms of the output measured over the interval $[t_k, t_{k+1})$.

\emph{Mapping for the differentiable part $z_k^\diff$:} In order to define $z_k^\diff \in \R^{r_k}$, where $r_k = \rk O_k^\diff$, we first introduce the function $\mathbf{z}_k^\diff:(t_k,t_{k+1}) \rightarrow \R^{r_k}, t\mapsto {Z_k^\diff}^\top e(t)$, which represents the observable component of the subsystem $(E_k, A_k, C_k)$ that can be recovered from the smooth output measurements $y^e$ over the interval $(t_k,t_{k+1})$. It follows (cf.\ \cite[Lem.~17]{TanwTren17a}) that the evolution of $\mathbf{z}_k^\diff$ is governed by an observable ODE
\begin{equation}\label{eq:zsys}
\begin{aligned}
 \dot{\mathbf{z}}_k^\diff &= S_k^\diff \mathbf{z}_k^\diff, \\
  y^e &= R_k^\diff \mathbf{z}_k^\diff, 
\end{aligned}
\end{equation}
where $S_k^\diff:={Z_k^\diff}^\top A_k^\diff Z_k^\diff$ and $R_k^\diff:=C_k^\diff Z_k^\diff$. Because of the observability of the pair $(S_k^\diff, R_k^\diff)$ in \eqref{eq:zsys}, there exists a (linear) operator $\fO_{(t_k,t_{k+1})}^\diff$ such that
\[
\mathbf{z}_k^\diff = \fO_{(t_k,t_{k+1})}^\diff (y^e_{(t_k,t_{k+1})})
\]
and we set
\[
z_k^\diff = \mathbf{z}_k^\diff(t_k^+).
\]
Note that for an estimation $\widehat{\mathbf{z}}_k^\diff$ of $\mathbf{z}_k^\diff$ obtained by a standard observer (e.g.\ the Luenberger observer) the evaluation at the beginning of the observation interval is not meaningful (because this value is not affected by the output injection). However, a good estimate can easily be obtained by propagating back the final estimate with the known homogeneous system dynamics; i.e.\
\begin{equation}\label{eq:hatzkdiff}
   \widehat{z}_k^\diff = \me^{-S_k^\diff \tau_k} \widehat{\mathbf{z}}_k^\diff(t_{k+1}^-).
\end{equation}

\emph{Mapping for the impulsive part $z_k^\imp$:}
The impulsive part of the output at switching time $t_k$ can be represented as
\[
   y^e[t_k] = \sum_{j=0}^{n-2} \eta_k^j \delta_{t_k}^{(j)},
\]
where due to Lemma~\ref{lem:impulses} the coefficients $\eta_k^j$ satisfy the relation $\boldsymbol{\eta}_k = -O_{k}^{\imp} e(t_k^-)$, with $\boldsymbol{\eta}_k := ({\eta_k^0}/\cdots/{\eta_k^{n-2}})\in\R^{(n-1) \dy}$.
We chose a matrix $U_k^\imp$ such that
\begin{equation}\label{eq:Ukimp}
   -{O_{k}^{\imp}}^\top U_k^\imp = Z_k^\imp,
\end{equation}
 then
\begin{equation}\label{eq:zkimp}
  z_k^\imp\!=\!{Z_k^\imp}^\top\! e(t_k^-)\!=\!-{U_k^\imp}^\top\! O_{k}^{\imp}  e(t_k^-)\!=\!{U_k^\imp}^\top\! \boldsymbol{\eta}_k.
\end{equation}

\subsection{Observable component over an interval}\label{sec:intObsComp}
For $q-1 \ge k \ge p$, the $[t_k,t_q)$-unobservable subspace \eqref{eq:Nkq} can be computed recursively as follows
\begin{subequations}\label{eq:DAESeqUnobs}
\begin{align}
\cN_{q-1}^{q} & = \cW_{q-1} \label{eq:DAESeqUnobsa}\\
\cN_{k}^q &= \cW_k \cap \Pi_{k}^{-1}\me^{-A_k^\diff \tau_k}\cN_{k+1}^q, \quad k\leq q-2. \label{eq:DAESeqUnobsb}
\end{align}
\end{subequations}

The objective is to compute the observable part $\mu_k = {M^q_k}^\top e (t_k^-)$ in \eqref{eq:MkNk} recursively for $k=q-1,q-2,\ldots,p$. We choose $\mu_{q-1} = z_{k-1}$.
By construction, we know that
\begin{align*}
\im M^q_k &= \cM^q_k = (\cN^q_k)^\bot = \left( \cW_k \cap \Pi_k^{-1}(\me^{-A_k^\diff \tau_k}\cN_{k+1}^{q}) \right)^\perp \\
& = \cW_k^\perp + \Pi_k^\top (\me^{-A_k^\diff \tau_k} \cN_{k+1}^{q})^\perp, \quad k\leq q-2.
\end{align*}
Recalling that $\im Z_k = (\cW_k)^\perp$, and introducing the matrix $\Theta^q_k$, for $k=p,p+1,\ldots,q-2$, such that
\begin{equation}\label{eq:Thetaqk_def}
\im \Theta^q_k = (\me^{-A_k^\diff \tau_k} \cN_{k+1}^{q})^\perp
\end{equation}
we obtain
\begin{equation}\label{eq:Mqk_def}
\im M^q_k = \im \left[Z_k, \ \Pi_k^\top \Theta^q_{k} \right].
\end{equation}
Hence there exists a matrix $U^q_k$ such that
\begin{equation}\label{eq:Uqk_def}
M^q_k = \left[Z_k, \ \Pi_k^\top \Theta^q_{k} \right] U^q_{k}.
\end{equation}
Noting that
\begin{align*}
\Pi_k e(t_k^-) &= e(t_k^+) = \me^{-A_k^\diff \tau_k} e(t_{k+1}^-)\\
& = \me^{-A_k^\diff \tau_k} \left( M^q_{k+1} \mu_{k+1} + N^q_{k+1} \nu_{k+1} \right)
\end{align*}
and multiplication on both sides from left by ${\Theta^q_{k}}^\top$ gives
\begin{multline*}
{\Theta^q_{k}}^\top \Pi_k e(t_k^-) = {\Theta^q_{k}}^\top \me^{-A_k^\diff \tau_k} M^q_{k+1} \mu_{k+1} \\
+ \underbrace{{\Theta^q_{k}}^\top \me^{-A_k^\diff \tau_k} N^q_{k+1}}_{=0} \nu_{k+1}.
\end{multline*}
This allows us to compute $\mu_k$, $k=q-2, q-3,\ldots,p$, as follows:
\begin{multline}\label{eq:muk_recursion}
\mu_{k} = {M^q_{k}}^\top e(t_k^-) = {U^q_{k}}^\top \begin{bmatrix} Z_k^\top \\ {\Theta^q_{k}}^\top \Pi_k \end{bmatrix} e(t_k^-) \\
 \!\!= {U^q_{k}}^\top\!\! \begin {pmatrix} z_k \\ {\Theta^q_{k}}^\top \me^{-A_k^\diff \tau_k} M^q_{k+1} \mu_{k+1} \end{pmatrix}
 =: \cF^q_k \begin{pmatrix} z_k \\ \mu_{k+1}\end{pmatrix}\!.
\end{multline}

%

\section*{References}

\begin{thebibliography}{0}
\expandafter\ifx\csname natexlab\endcsname\relax\def\natexlab#1{#1}\fi
\expandafter\ifx\csname url\endcsname\relax
  \def\url#1{\texttt{#1}}\fi
\expandafter\ifx\csname urlprefix\endcsname\relax\def\urlprefix{URL }\fi

\end{thebibliography}


\begin{thebibliography}{22}
\expandafter\ifx\csname natexlab\endcsname\relax\def\natexlab#1{#1}\fi
\expandafter\ifx\csname url\endcsname\relax
  \def\url#1{\texttt{#1}}\fi
\expandafter\ifx\csname urlprefix\endcsname\relax\def\urlprefix{URL }\fi

\bibitem[{Astolfi and Praly(2003)}]{AstoPral03}
Astolfi, A., Praly, L., 2003. Global complete observability and output-to-state
  stability imply the existence of a globally convergent observer. In: Proc.
  42nd~{IEEE} Conf. Decis. Control, Hawaii, USA. pp. 1562--1567.

\bibitem[{Berger et~al.(2012)Berger, Ilchmann, and Trenn}]{BergIlch12a}
Berger, T., Ilchmann, A., Trenn, S., 2012. The quasi-{W}eierstra{\ss} form for
  regular matrix pencils. Linear Algebra Appl. 436~(10), 4052--4069.

\bibitem[{Berger and Reis(2017)}]{BergReis17c}
Berger, T., Reis, T., 2017. Observers and dynamic controllers for linear
  differential-algebraic systems. {SIAM} J. Control Optim. 55~(6), 3564--3591.

\bibitem[{Berger et~al.(2017)Berger, Reis, and Trenn}]{BergReis17a}
Berger, T., Reis, T., Trenn, S., 2017. Observability of linear
  differential-algebraic systems: A survey. In: Ilchmann, A., Reis, T. (Eds.),
  Surveys in {D}ifferential-{A}lgebraic {E}quations {IV}.
  Differential-Algebraic Equations Forum. Springer-Verlag, Berlin-Heidelberg,
  pp. 161--219.

\bibitem[{De~Santis et~al.(2009)De~Santis, Di~Benedetto, and Pola}]{DeSaDiBe09}
De~Santis, E., Di~Benedetto, M.~D., Pola, G., 2009. A structural approach to
  detectability for a class of hybrid systems. Automatica 45~(5), 1202--1206.

\bibitem[{Liberzon and Trenn(2009)}]{LibeTren09}
Liberzon, D., Trenn, S., December 2009. On stability of linear switched
  differential algebraic equations. In: Proc. {IEEE} 48th Conf. on Decision and
  Control. pp. 2156--2161.

\bibitem[{Liberzon and Trenn(2012)}]{LibeTren12}
Liberzon, D., Trenn, S., May 2012. Switched nonlinear differential algebraic
  equations: Solution theory, {L}yapunov functions, and stability. Automatica
  48~(5), 954--963.
  
\bibitem[Mancilla-Aguilar and Garc\`ia(2018)]{MancGarc18}
Mancilla-Aguilar, J.L., Garc\`ia, R.A., 2018. Uniform Asymptotic Stability of 
Switched Systems via detectability of reduced control systems. In: Proc. 2018
American Control Conference. pp. 4552--4557.

\bibitem[{M{\"u}ller and Liberzon(2012)}]{MullLibe12}
M{\"u}ller, M., Liberzon, D., 2012. Input/output-to-state stability and
  state-norm estimators for switched nonlinear systems. Automatica 48~(9),
  2029--2039.

\bibitem[{Schwartz(1950, 1951)}]{Schw50}
Schwartz, L., 1950, 1951. Th\'{e}orie des Distributions I,II. No. IX,X in
  Publications de l'institut de math{\'e}matique de l'Universite de Strasbourg.
  Hermann, Paris.

\bibitem[{Shim et~al.(2012)Shim, Tanwani, and Ping}]{ShimTanw12}
Shim, H., Tanwani, A., Ping, Z., December 2012. Back-and-forth operation of
  state observers and norm estimation of estimation error. In: Proc.
  51st~{IEEE} Conf. Decis. Control, Maui, USA. pp. 3221--3226.

\bibitem[{Sontag and Wang(1997)}]{SontWang97}
Sontag, E., Wang, Y., 1997. Output-to-state stability and detectability of
  nonlinear systems. Syst. Control Lett. 29~(5), 279--290.

\bibitem[{Sontag(1998)}]{Sont98a}
Sontag, E.~D., 1998. Mathematical Control Theory: Deterministic Finite
  Dimensional Systems, 2nd Edition. Springer-Verlag, New York.

\bibitem[{Tanwani et~al.(2013)Tanwani, Shim, and Liberzon}]{TanwShim13}
Tanwani, A., Shim, H., Liberzon, D., 2013. Observability for switched linear
  systems: Characterization and observer design. {IEEE} Trans. Autom. Control
  58~(4), 891--904.

\bibitem[{Tanwani and Trenn(2010)}]{TanwTren10}
Tanwani, A., Trenn, S., 2010. On observability of switched
  differential-algebraic equations. In: Proc. 49th~{IEEE} Conf. Decis. Control,
  Atlanta, USA. pp. 5656--5661.

\bibitem[{Tanwani and Trenn(2012)}]{TanwTren12}
Tanwani, A., Trenn, S., 2012. Observability of switched differential-algebraic
  equations for general switching signals. In: Proc. 51st~{IEEE} Conf. Decis.
  Control, Maui, USA. pp. 2648--2653.

\bibitem[{Tanwani and Trenn(2013)}]{TanwTren13}
Tanwani, A., Trenn, S., 2013. An observer for switched differential-algebraic
  equations based on geometric characterization of observability. In: Proc.
  52nd~{IEEE} Conf. Decis. Control, Florence, Italy. pp. 5981--5986.

\bibitem[{Tanwani and Trenn(2015)}]{TanwTren15}
Tanwani, A., Trenn, S., 2015. On detectability of switched linear
  differential-algebraic equations. In: Proc. 54th~{IEEE} Conf. Decis. Control,
  Osaka, Japan. pp. 2957--2962.

\bibitem[{Tanwani and Trenn(2017{\natexlab{a}})}]{TanwTren17a}
Tanwani, A., Trenn, S., 2017{\natexlab{a}}. Determinability and state
  estimation for switched differential–algebraic equations. Automatica 76,
  17--31.

\bibitem[{Tanwani and Trenn(2017{\natexlab{b}})}]{TanwTren17b}
Tanwani, A., Trenn, S., 2017{\natexlab{b}}. Observer design for detectable
  switched differential-algebraic equations. In: Proceedings of the 20th IFAC
  World Congress. pp. 2953--2958.

\bibitem[{Trenn(2009)}]{Tren09d}
Trenn, S., 2009. Distributional differential algebraic equations. Ph.D. thesis,
  Institut f{\"u}r Mathematik, Technische Universit{\"a}t Ilmenau,
  Universit{\"a}tsverlag Ilmenau, Germany.
\newline\urlprefix\url{http://www.db-thueringen.de/servlets/DocumentServlet?id=13581}

\bibitem[{Trenn(2012)}]{Tren12}
Trenn, S., 2012. Switched differential algebraic equations. In: Vasca, F.,
  Iannelli, L. (Eds.), Dynamics and Control of Switched Electronic Systems -
  Advanced Perspectives for Modeling, Simulation and Control of Power
  Converters. Springer-Verlag, London, Ch.~6, pp. 189--216.
  
\bibitem[{Trentelman et~al.(2001)Trentelman, Stoorvogel, and
  Hautus}]{TrenStoo01}
Trentelman, H.~L., Stoorvogel, A.~A., Hautus, M. L.~J., 2001. Control Theory
  for Linear Systems. Communications and Control Engineering. Springer-Verlag,
  London.
  

\bibitem[{Wong(1974)}]{Wong74}
Wong, K.-T., 1974. The eigenvalue problem {$\lambda Tx + Sx $}. J. Diff. Eqns.
  16, 270--280.

\end{thebibliography}

\end{document}